\documentclass[12pt]{article}
\usepackage[english]{babel}
\usepackage[a4paper, margin=1in]{geometry}
\usepackage{natbib}
\usepackage{setspace}
\usepackage{amsmath, amssymb, bbm, bm}
\usepackage{graphicx}
\usepackage[colorlinks=true, allcolors=blue]{hyperref}
\usepackage[dvipsnames]{xcolor}
\usepackage{comment}

\title{Some benefits of standardisation for conditional extremes}
\author{Christian Rohrbeck$^{1}$\footnote{Corresponding author, cr777@bath.ac.uk}~~and Jonathan A. Tawn$^2$\\
\normalsize{$^1$Department of Mathematical Sciences,
University of Bath, BA2 7AY, U.K.}\\
\normalsize{$^2$Department of Mathematics and Statistics, Lancaster University, LA1 4YF, U.K.}}
\date{\vspace{-0.5cm}}

\begin{document}

\maketitle

\begin{abstract}
A key aspect where extreme values methods differ from standard statistical models is through having asymptotic theory to provide a theoretical justification for the nature of the models used for extrapolation. In multivariate extremes many different asymptotic theories have been proposed, partly as a consequence of the lack of ordering property with vector random variables. One class of multivariate models, based on conditional limit theory as one variable becomes extreme, developed by \citet{HT2004}, has developed wide practical usage. The underpinning value of this approach has been supported by further theoretical characterisations of the limiting relationships by \cite{HeffernanResnick2007} and \cite{Resnick2014}. However, \cite{Drees2017} provide a number of counterexamples to their results. This paper studies these counterexamples in the \citet{Keef2013a} framework, which involves marginal standardisation to a common exponentially decaying tailed marginal distribution. Our calculations show that some of the issues identified by \cite{Drees2017} can be addressed in this way.
\end{abstract}
\onehalfspacing

\section{Introduction}

Multivariate extreme value problems are important across a range of subject domains, such as sea levels \citep{Coles1994}, air pollution \citep {HT2004}, rainfall \citep{Davison2012} and river flow \citep{Engelke2020}, which all feature in influential discussion papers.
The typical formulation is to have $n$ independent and identically distributed replicate observations $(\mathbf{x}_1,\ldots ,\mathbf{x}_n)$, from a $d$-dimensional vector random variable $\mathbf{X}$ with unknown joint distribution $F_{\mathbf{X}}$. Here the aim is to estimate $\Pr(\mathbf{X}\in A)$ where $A\subset \mathbb{R}^d$, such that all the elements in $A$ are in the upper tail of at least one of the marginal distributions of $\mathbf{X}$, with the formulation of $A$ depending on the characteristics of the problem of interest. The typical approach to make such inference is to estimate both the marginal distributions and dependence structure (copula) with a focus on their behaviour in their upper extremes. Univariate extreme value methods are well established \citep{Coles2001, DavisonSmith1990}, with multivariate dependence modelling being the key challenge. 

In the bivariate case, that we will focus on for variables $(X, Y)$, there are two distinct types of extremal dependence, which are easiest explained via the coefficient of asymptotic dependence $\chi= \lim_{p\,\uparrow\, 1}\Pr\left[F_Y(Y)>p\mid F_X(X)>p\right]$,
where $F_X$ and $F_Y$ are the marginal distributions of $X$ and $Y$ respectively. Having $\chi>0$ coincides with asymptotic dependence between $X$ and $Y$, a situation in which both variables can take their largest values simultaneously; while when $\chi=0$, termed asymptotic independence, such limiting dependence is impossible. Many models for bivariate extremes are only suitable in one of these situations: multivariate max-stable distributions \citep{Gudendorf2012} and multivariate generalised Pareto distributions \citep{Kiriliouk2019} only allow $\chi>0$ or independent variables. Therefore, distinguishing between these cases, or having a model that incorporates both in a flexible way, can play a crucial role in model selection.

One class of multivariate extreme models, based on conditional limit theory as one variable becomes extreme, developed by \citet{HT2004}, has developed  wide practical usage, with applications linked to widespread river flooding \citep{Keef2013b}, time series dependence in heatwaves \citep{WinterTawn2017}, spatial air temperature extremes \citep{WadsworthTawn2022}, spatio-temporal sea-surface temperatures \citep{Simpson2021}, offshore metocean environmental design contours \citep{Ewans2014}, coastal flooding \citep{Gouldby2017}, 
food chemicals \citep{Paulo2006}, and laboratory trials \citep{Southworth2012}.

This \citet{HT2004} class of models has considerable flexibility as it covers both asymptotic dependence and asymptotic independence classes. Furthermore, in the multivariate case it allows for different extremal dependence classes between separate subsets of the variables, unlike models such as \cite{Wadsworth2017} and \cite{HuserWadsworth2019}. Since its initial presentation, the model proposed by \citet{HT2004} has been extended by \citet{Keef2013a} to its current most widely adopted form. Specifically, 
for $(X,Y)$ marginally transformed to have Laplace marginals, denoted $(X_L,Y_L)$, it is assumed that there exists values $(\alpha_{\mid X},\beta_{\mid X})\in[-1,1]\times(-\infty,1)$ such that for $x>0$ and $z\in \mathbb{R}$
\begin{equation}
\Pr\left\{\frac{Y_L-\alpha_{\mid X}X_L}{X_L^{\beta_{\mid X}}} \leq z, X_L-t>x \mid X_L>t\right\} \to G_{\mid X}(z)\exp(-x) \qquad \mathrm{as}~t\to\infty,
\label{eq:HT2004NEW}
\end{equation} 
where $G_{\mid X}$ is the distribution function of a non-degenerate random variable with values in $[-\infty,\infty)$; the restriction that $G_{\mid X}$ places no mass at $\{\infty\}$ ensures that $\alpha_{\mid X}$ is uniquely identifiable. This relation gives that the normalised $Y_L$ is conditionally independent of $X_L$ in the limit. To characterise the full joint tail of $(X_L,Y_L)$ in addition to limit~\eqref{eq:HT2004NEW} we also need the equivalent relationship for the reverse conditional distribution of $X_L$ given $Y_L$ is large.

Despite the strong applied value of the conditional modelling framework, some concerns about the broader theoretical restrictions of the limiting assumptions exist. Attempts to formalise the method and weaken some of these assumptions include \cite{HeffernanResnick2007} and \cite{Resnick2014}. These results involve the both random and non-random norming, which in  the context of limit~\eqref{eq:HT2004NEW} corresponds to norming $Y_L$ by a function of $X_L$ or by a function of $t$ respectively. However, \cite{Drees2017} provided a number of counterexamples of their results.

This paper explores these counterexamples to see if they undermine any of the asymptotic justifications for the statistical methods stemming from the \citet{HT2004} framework and their practical adoption. There is a critical difference between the framework studied in \cite{HeffernanResnick2007}, \cite{Resnick2014}, \cite{Drees2017} from the \citet{HT2004} framework, specifically, that latter requires an initial marginal standardisation, so that after transformation of $(X,Y)$ they are assumed to have identical marginal distributions before studying the conditional extremes behaviour. This transformation was taken to be Gumbel in \citet{HT2004} and Laplace (as above) in \citet{Keef2013a}. Such standardisation of variables to common margins is quite usual in the study of dependence structure, e.g.,  \citet{Nelsen1999} and \citet{Beirlant2004}, as this makes  relationships more easy to model through linearity, with exponential margins being particularly desirable for this, as shown by
\citet{Papastathopoulos2017}. Our intuition is that having marginal variables on completely different marginal tail behaviours (explicitly different shape parameters/tail indices) imposes a major restriction on a conditional approach using affine transformations, such as in the norming of $Y_L$ in limit~\eqref{eq:HT2004NEW}. 
Furthermore, we believe that there are advantages of random over non-random norming, as it is the exact value of the conditioning variable that affects the response variable, which in bivariate cases leads to simpler limits and in multivariate cases can lead to conditional independence which is vital for the study of graphical structures 
\citep{Engelke2020}.

The paper is structured as follows: In Section~\ref{sec:background} we present the background theory of the different conditional representations. In Section~\ref{sec:Examples} we cover each of the counterexamples given by \cite{Drees2017}, with simulations to help interpretation, and state which features of \cite{Das2011} and \cite{Resnick2014} they show are not appropriate. In each case we illustrate how the problems are overcome through an initial standardisation of the marginal distributions. We then discuss the practical implications of our calculations in Section~\ref{sec:Implications} and conclude with a discussion in Section~\ref{sec:Discussion}. Some technical details of the calculations for the examples are given in the Appendix.

\section{Background Theory}
\label{sec:background}

\subsection{Multivariate and Conditional Extremes}

%The first models studied in multivariate extreme value theory considered the limiting behaviour of the joint distribution of componentwise maxima of vector random variables. 
For notational simplicity, we focus on the bivariate case with $(X,Y)$, where $X$ and $Y$ are continuous random variables. Classical multivariate extreme value models assume that the marginal distributions $F_X$ and $F_Y$ of $(X,Y)$ belong to the domain of attraction of some extreme value distribution: $F_X$ is in the domain of attraction of an extreme value distribution if there exist functions $p_X:\mathbb{R}\to\mathbb{R}_+$ and $q_X:\mathbb{R}\to\mathbb{R}$ such that
\begin{equation}
F_X^t\left\{p_X(t)x + q_X(t)\right\} \to \exp\left\{-(1+\gamma_X x)^{-1/\gamma_X}\right\}\quad \mathrm{as}~t\to\infty
\label{eq:Tail}
\end{equation}
for some $\gamma_X\in\mathbb{R}$ and all $x\in E^{(\gamma_X)}:=\{x \in\ \mathbb{R} \mid 1 + \gamma_X x >0\}$. Multivariate extreme value distributions then arise as the limiting joint distribution of the componentwise maxima of independent and identically distributed random variables $(X_i,Y_i)$, for $i=1,\ldots,t$, with joint distribution function $F_{X,Y}$ and marginal distribution functions 
$X_i\sim F_X$ and $Y_i\sim F_Y$. Specifically, it as assumed that there exist functions $p_X, q_X$ as in limit~\eqref{eq:Tail}, and similarly $p_Y, q_Y$, such that 
\begin{align*}
\Pr &\left(\frac{\max_{i=1,\ldots,t}X_i-q_X(t)}{p_X(t)} \leq x, \frac{\max_{i=1,\ldots,t} Y_i-q_Y(t)}{p_Y(t)}\leq y\right)\\
= & [F_{X,Y}(p_X(t)x+q_X(t), p_Y(t)y+q_Y(t))]^t\to H(x,y) \qquad \mbox{as}~t\to\infty, 
\end{align*}
where $H$ is a bivariate distribution function with non-degenerate marginal distributions, given by limit form~\eqref{eq:Tail}, with tail indices of $\gamma_X$ and 
$\gamma_Y$ respectively,  and with a copula possessing a specific max-stable property which, amongst other features, excludes the possibility of negative dependence, see \citet{ColesHeffernanTawn1999} and \citet{Beirlant2004}.

\citet{HT2004} propose examining the dependence in the tail of ($X,Y$) by first standardising the marginals via the probability integral transformation to have Gumbel distributions, denoted $(X_G,Y_G)$, with $\Pr(X_G\le x)=\Pr(Y_G\le x)=\exp\{-\exp(-x)\}$ for $x\in\mathbb{R}$, and considering the conditional distribution of $Y_G\mid (X_G=t)$ as $t\to\infty$. The assumption underlying their approach is that there exist normalising functions $\tilde{a}_{\mid X}(y):\mathbb{R}_+\to\mathbb{R}$ and $\tilde{b}_{\mid  X}(y):\mathbb{R}_+\to\mathbb{R}_+$ such that
\begin{equation}
\Pr\left\{ \frac{Y_G-\tilde{a}_{\mid  X}(X_G)}{\tilde{b}_{\mid X}(X_G)} \leq z \mid X_G=t\right\} \to \tilde{G}_{\mid X}(z) \qquad \mathrm{as}~t\to\infty,
\label{eq:HT2004}
\end{equation}
where the limit distribution $\tilde{G}_{\mid X}$ is non-degenerate. To ensure that $\tilde{a}_{\mid X}$, $\tilde{b}_{\mid X}$ and $\tilde{G}_{\mid X}$ are well-defined, we require $\lim_{z\to\infty}\tilde{G}_{\mid X}(z)=1$, i.e., $\tilde{G}_{\mid X}$ has no mass at $+\infty$, and $\tilde{b}_{\mid X}(x)/x\rightarrow 0$ as $x\rightarrow \infty$ \citep{Keef2013a}. \citet{HT2004} find that (up to type) the functions $\tilde{a}_{\mid X}$ and $\tilde{b}_{\mid X}$ in~\eqref{eq:HT2004} have a common parametric form for all copulas described by \cite{Joe1997} and \cite{Nelsen1999}.

\cite{HeffernanResnick2007} modify and extend the original framework by \cite{HT2004}. Their main modification is to replace the condition $X=t$ in \eqref{eq:HT2004} by $X>t$, i.e., they analyse
\begin{equation}
\Pr\left\{ \frac{Y-a_{\mid  X}(X)}{b_{\mid X}(X)} \leq z \mid X>t\right\} \to G_{\mid X}(z) \qquad \mathrm{as}~t\to\infty.
\label{eq:HR20071}
\end{equation}
This is the most widely applied and considered conditional extreme value model framework and we will use it in the remainder of the paper. \cite{HeffernanResnick2007} further drop the assumption that $X$ and $Y$ have Gumbel margins and provide theoretical results subject to $F_X$ lying in the domain of attraction of some extreme value distribution. 

\cite{Keef2013a} focus on $X$ and $Y$ having standard Laplace margins, i.e.,
\[
\Pr(X<x) = 
\begin{cases}
\frac{1}{2} \exp(x)&\mbox{if}~ x\leq 0,\\
1 - \frac{1}{2} \exp(-x)&\mbox{if}~ x>0.
\end{cases}
\]
Under these conditions, the functions in \eqref{eq:HR20071} are of the form $a_{\mid X}(x) = \alpha x$ and $b_{\mid X}(x) = x^\beta$, with $(\alpha,\beta)\in[-1,1]\times(-\infty,1)$, in all of the standard copulas studied by \cite{HT2004}. When $X$ and $Y$ are positively associated, standardisation of $X$ and $Y$ to Laplace margins gives the same limiting behaviour as when the variables were transformed to Gumbel margins. However, the limiting behaviours differ when $X$ and $Y$ are negatively associated, with the symmetry of the Laplace margins giving a simpler form. Estimation of the parameters $(\alpha,\beta)$ and the distribution function $G_{\mid X}$ for Laplace margins are considered in \cite{Keef2013a, Keef2013b}. 

\subsection{Linking Multivariate and Conditional Extremes Models}

One early question asked about the conditional extremes model by \citet{HT2004} concerned its link to established multivariate extreme value models, e.g., \citet{Husler1989} or \citet{Tawn1990}. This motivated \cite{HeffernanResnick2007} to define the class of conditional extreme value models (CEVM) which does not require the margins to be standardised to common margins. The limit distribution of $Y\mid (X>t)$ as $t\to\infty$ lies in the CEVM class if
\begin{enumerate}
\item The distribution function $F_X$ of $X$ is in a domain of attraction of an extreme value distribution with parameter $\gamma_X\in\mathbb{R}$.

\item There exist normalising functions $c,f: \mathbb{R} \to \mathbb{R}$ and $d,g: \mathbb{R} \to \mathbb{R}_+$, such that 
\begin{equation}
t\,\Pr\left\{\frac{X-c(t)}{d(t)}>x, \frac{Y-f(t)}{g(t)}\leq y\right\} \to \mu_{Y\mid X>}\left((x,\infty]\times[-\infty,y]\right)
\quad\mathrm{as}~t\to\infty,
\label{eq:HR2007}
\end{equation}
where $\mu_{Y\mid X>}\left((x,\infty] \times [-\infty,y]\right)$ is a non-degenerate distribution function in $y$, and $\mu_{Y\mid X>}\left((x,\infty] \times [-\infty,y]\right)<\infty$.
\end{enumerate}
\cite{Resnick2014} add the condition $\mu_{Y\mid X>}((x,\infty]\times\{\infty\})=0$ to ensure uniqueness of the limit measure. However, Example 2.3 in \cite{Drees2017} show that this condition has to be strengthened to 
\[
\lim_{y\rightarrow \infty}\mu_{Y\mid X>}((x,\infty]\times [y,\infty) )
=
\lim_{y\rightarrow \infty}\mu_{Y\mid X>}((x,\infty]\times (-\infty,-y])
=0,
\]
i.e., the limit measure cannot put any mass for $Y$ at $\{-\infty\}$ or $\{+\infty\}$.

Returning to the link between the CEVM framework and multivariate extreme value models assume that $X$ and $Y$ belong to the domain of attraction of some extreme value distribution, with parameters $\gamma_X$ and $\gamma_Y$ respectively. Theorem 2.1 in \cite{Das2011} states that ($X,Y$) lies in the domain of attraction of a multivariate extreme value distribution if $Y\mid (X>t)$ and $X\mid (Y>t)$ both lie in the CEVM class by \cite{HeffernanResnick2007}. Example 4.4 in \cite{Drees2017} illustrates, however, that this result is not true, unless the normalisations of $Y\mid(X>t)$ and $X\mid(Y>t)$ in \eqref{eq:HR2007} are identical. 

Now consider the case that $(X,Y)$ lies in the domain of attraction of a multivariate extreme value distribution. Suppose that (i) $X$ and $Y$ are asymptotically dependent and (ii) $\gamma_X,\gamma_Y\leq 0$. \cite{Drees2017} show that these conditions are sufficient for the limits of $Y\mid (X>t)$ and $X\mid (Y>t)$ to lie in the class of conditional extreme value models. The restriction $\gamma_X,\gamma_Y\leq 0$ is required, as demonstrated by Example 4.2 in \cite{Drees2017}. Note, the conditions (i) and (ii) are not necessary conditions for $Y\mid (X>t)$ and $X\mid (Y>t)$ to lie in the CEVM class; see Section 5 in \cite{HeffernanResnick2007} for the case of $X$ and $Y$ being asymptotically independent.

\subsection{Standardisation of Marginals in the CEVM Class}

\cite{HeffernanResnick2007}  examined how the standardisation of $X$ to a standard Pareto distributed random variable $X_P$, but leaving $Y$ unchanged, affected the limiting measure $\mu_{Y\mid X>}$ in~\eqref{eq:HR2007}. They show that the limiting behaviour of $(X_P,Y)$ satisfies
\[
t\,\Pr\left\{\frac{X_P}{t}>x, \frac{Y-f(t)}{g(t)}\leq y\right\} \to 
\begin{cases}
\mu_{Y\mid X>}\left(\left(\dfrac{x^{\gamma_X}-1}{\gamma_X},\infty\right]\times[-\infty,y]\right) & \mbox{if}~\gamma_X\neq 0,\\
\mu_{Y\mid X>}\left((\log x,\infty]\times[-\infty,y]\right) & \mbox{if}~\gamma_X=0,
\end{cases}
\]
where the measure $\mu_{Y\mid X>}$ corresponds to that in \eqref{eq:HR2007}, and where $c(t)=0$ and $d(t)=t$ due to the standardisation to $X_P$.

The standardisation of $Y$ is more challenging than that of $X$, because the CEVM~\eqref{eq:HR2007} does not require $F_Y$ to be in a domain of attraction of an extreme value distribution, unlike the
\citet{HT2004} and \cite{Keef2013a} formulations of limit~\eqref{eq:HT2004}. \citet{HeffernanResnick2007} and \cite{Das2011} consider the task of finding a monotone and unbounded function $h:\mathbb{R}\to\mathbb{R}_+$ such that 
\begin{equation}
t\,\Pr\left\{\frac{X_P}{t}>x, \frac{h(Y)}{t}\leq y \right\} \to \tilde{\mu}_{Y\mid X>}\left((x,\infty]\times[-\infty,y]\right)
\quad\mathrm{as}~t\to\infty,
\label{eq:Standardization}
\end{equation}
where $\tilde{\mu}_{Y\mid X>}$ is finite and non-degenerate, $\tilde\mu_{Y\mid X>}((x,\infty]\times\{\infty\})=0$ and $\tilde\mu_{Y\mid X>}(\{\infty\}\times[-\infty,y])=0$ for all $x,y$. \cite{Das2011} argue that such a function $h$ exists if, and only if, $\mu_{Y\mid X>}$ is not a product measure. However, examples in \cite{Drees2017} Section 3 illustrate that neither implication is true, and the limit measures $\mu_{Y\mid X>}$ and $\tilde{\mu}_{Y\mid X>}$ may convey different information if they exist. \cite{Drees2017} further provide two sufficient sets of conditions on the functions $f$ and $g$ in expression~\eqref{eq:HR2007} for such a function $h$ to exist. 

%The first set of conditions is $g(t)\to\infty$ as $t\to\infty$, $f(t)=0$ and $\mu_{Y\mid X>}((x,\infty)\times(0,\infty))>0$ for all $x>0$, and the second set of conditions is $g(t)\to 0$ as $t\to\infty$, $f(t)=c_0$ for some $c_0\in\mathbb{R}$, and the range of $Y$ is bounded by $c_0$, range$(Y)\subset (-\infty,c_0)$ or range$(Y)\subset (c_0, \infty)$.

\section{Investigating \citet{Drees2017} Examples}
\label{sec:Examples}
\subsection{Strategy}
\label{sec:Strategy}

Most examples in \citet{Drees2017} work with the joint distribution of $(X_P,Y)$, where $\Pr(X_P>x) = 1-x^{-1}$ ($x>1$), i.e., the conditioning variable has standard Pareto distribution, and the distribution of $Y$ is given indirectly through the distributions of $X_P$ and $Y\mid X_P$. In the following, we consider the examples by \cite{Drees2017} highlighted in Section~\ref{sec:background} and examine the obtained limiting behaviour, in each case using their numbering of the examples. We work within the framework of \cite{Keef2013a}, therefore, we consider the limiting behaviour after the variables $X_P$ and $Y$ have been transformed to Laplace margins, denoted by $(X_L,Y_L)$. 

The standardisation of $X_P$ to Laplace margins yields the variable $X_L$ and the link between the values $x_L$ of $X_L$ and $x$ of $X_P$ is given by
\begin{equation}
\frac{1}{x} =
\begin{cases}
1-\frac{1}{2}\exp(x_L) & \mathrm{if}~x_L\leq 0,\\
\frac{1}{2}\exp(-x_L) & \mathrm{if}~x_L>0.
\end{cases}
\label{eq:LaplaceX}
\end{equation}
When transforming $Y$ to $Y_L$, we first derive the distribution function $F_Y$ and then derive the expression for the transformed value $y_L$ of $Y_L$ as $y_L = F_L^{-1} \left[ F_Y(y)\right]$, where $F_L^{-1}$ is the inverse distribution function of a Laplace random variable. 

\subsection{Example 2.3}
Let $B$ be a discrete random variable that is uniformly distributed on $\{0,1\}$ that is independent of~$X_P$ and define $Y = B + (1-B)(2-1/X_P)$. The variable $Y$ can take any value in the interval $[1,2)$, with the highest values occurring when $B=0$ and $X_P$ large, and its marginal distribution is given by $\Pr(Y=1)=1/2$ and $\Pr(Y<y)= y/2$ for $1<y\leq 2$. This example in \cite{Drees2017} showed that the condition $\mu_{Y\mid X}\left([x,\infty]\times\{\infty\}\right)=0$ in \cite{Resnick2014} is not sufficient to ensure uniqueness of the limit measure in expression~\eqref{eq:HR2007}, and that the stronger condition $\mu_{Y\mid X}\left([x,\infty]\times\{-\infty,\infty\}\right)=0$ is required. 

\begin{figure}
\centering
\includegraphics[width=0.4\textwidth]{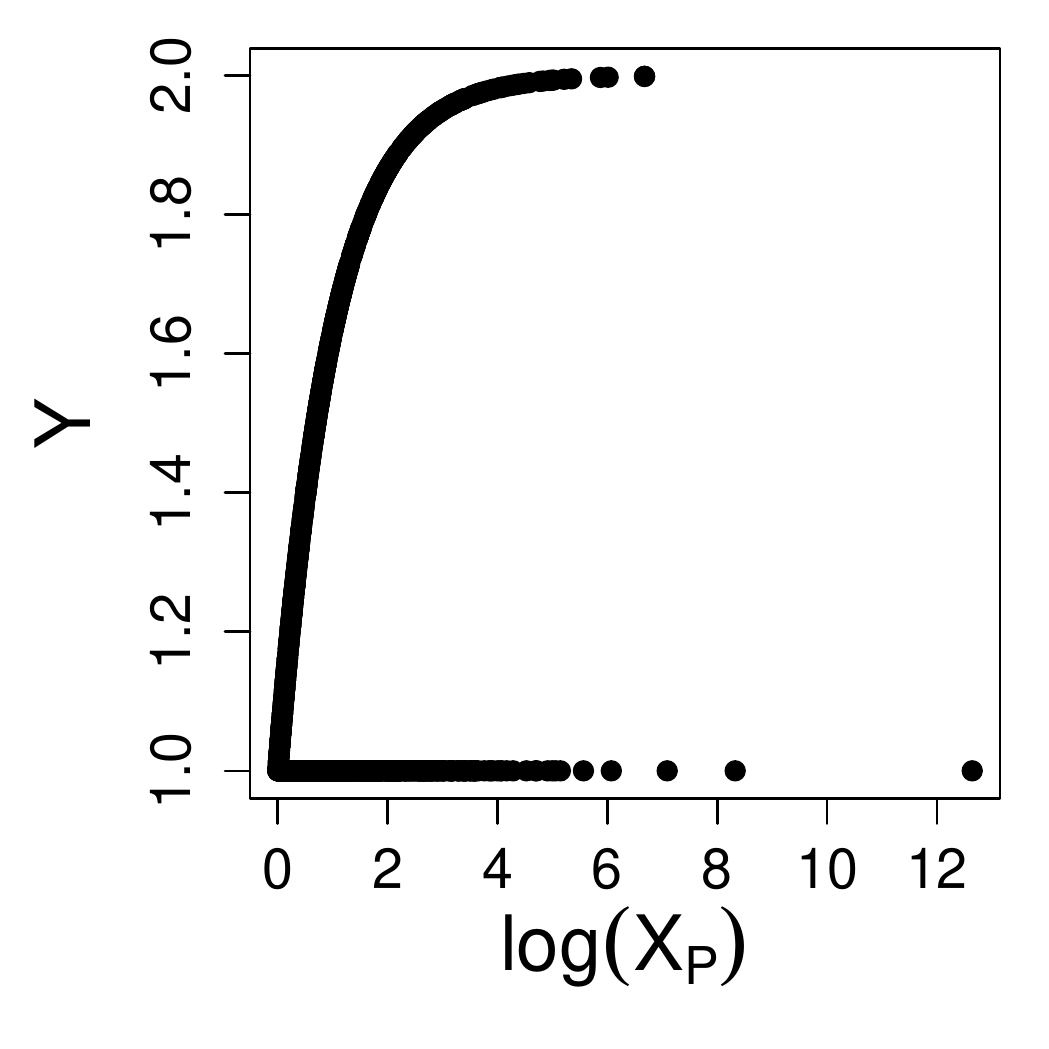}
\hspace{0.5cm}
\includegraphics[width=0.4\textwidth]{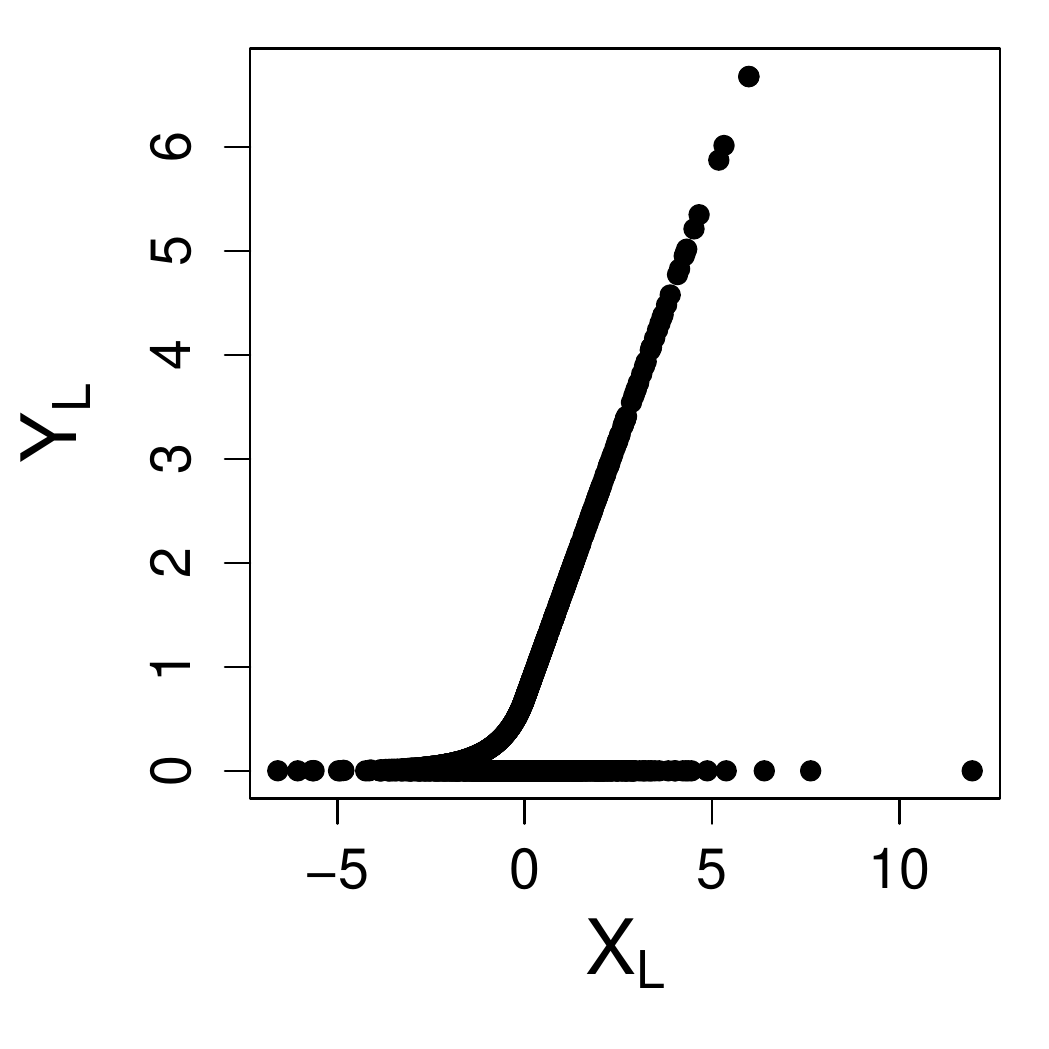}
\caption{Illustration of 2,000 samples for the framework in Example 2.3. The left panel shows the simulated observations $(x_p,y)$ on the original scale, while the right panel corresponds to the transformed samples $(x_L, y_L)$.}
\label{fig:Example23}
\end{figure}

As outlined in Section \ref{sec:Strategy}, we are interested in the limiting behaviour of the transformed variable $Y_L$ given that the transformed variable $X_L$ is large. Transformation of $Y$ to Laplace margins gives $Y = 2-\exp(-Y_L)$ for $B=0$, while $Y=1$ when $B=1$; this second case implies $Y_L=0$ irrespective of $X$ for $B=1$ and, thus, the lower tail of $Y_L$ is not Laplace distributed. Figure~\ref{fig:Example23} left panel shows that realised values of $Y$ are close to $y=1$ or $y=2$ for large values of $X_P$, while the transformed variable $Y_L$, shown in the right panel, has no upper bound. When $B=0$, substituting the realisations $x$ and $y$ by their transformed values $x_L$ and $y_L$, gives for $x_L>0$
\[
y = 2 - \frac{1}{x}
\quad\Leftrightarrow\quad 
2-\exp(-y_L) =~2 - \frac{1}{2}\exp(-x_L)
\quad\Leftrightarrow\quad 
y_L =~\log(2) + x_L. 
\]
This linear relationship between the values, for $x_L>0$ and $B=0$, is also visible in Figure \ref{fig:Example23} right panel, while, for $B=1$, we have $y_L=0$ for all possible values~$x_L$. 

From the calculations above, we conclude that the functions $a_{\mid X}(x) = x$ and $b_{\mid X}(x) = 1$ in expression~\eqref{eq:HR20071} give the limiting behaviour as
\[  
\Pr(Y_L - X_L \leq z \mid X_L > x_L) ~\to~ \frac{1}{2}\left(1 + \mathbf{I}\{\log 2 \leq z\}\right) = G_{\mid X}(z) \qquad \mathrm{as}~x_L \to \infty,
\]
where $\mathbf{I}$ denotes the indicator function, and $G_{\mid X}$ is a non-degenerate distribution function. The result $\lim_{z\to-\infty}G_{\mid X}(z)=0.5$ is due to the case $B=1$ which occurs with probability $0.5$. Other choices for $a_{\mid X}$ and $b_{\mid X}$ lead to a degenerate limiting distribution $G_{\mid X}$, contradicting the \cite{HT2004} assumption, or yield $\mu_{Y\mid X>}\left([x,\infty]\times\{\infty\}\right)=0$, violating the constraint $\lim_{z\to\infty} G_{\mid X}(z) =1$ by \cite{Keef2013a}.

In terms of the measure $\mu_{Y_L\mid X_L>}$ in \eqref{eq:HR2007}, we have for $x_L>0$
\[
\mu_{Y_L\mid X_L>}\left((x_L,\infty]\times [-\infty,y_L]\right) = \frac{1}{2} \left(1 + \mathbf{I}\{\log 2\leq y_L\}\right)\times \frac{1}{2} \exp(-x_L). 
\]
While this result is similar to the first limit found by \cite{Drees2017}, we do not require the additional constraint $\mu_{Y\mid X>}\left([x,\infty]\times\{-\infty\}\right)=0$, introduced by \cite{Drees2017}, to ensure a unique limiting behaviour, because we transformed the variables to common Laplace margins.

\subsection{Example 3.1}
\begin{figure}
\centering
\includegraphics[width=0.4\textwidth]{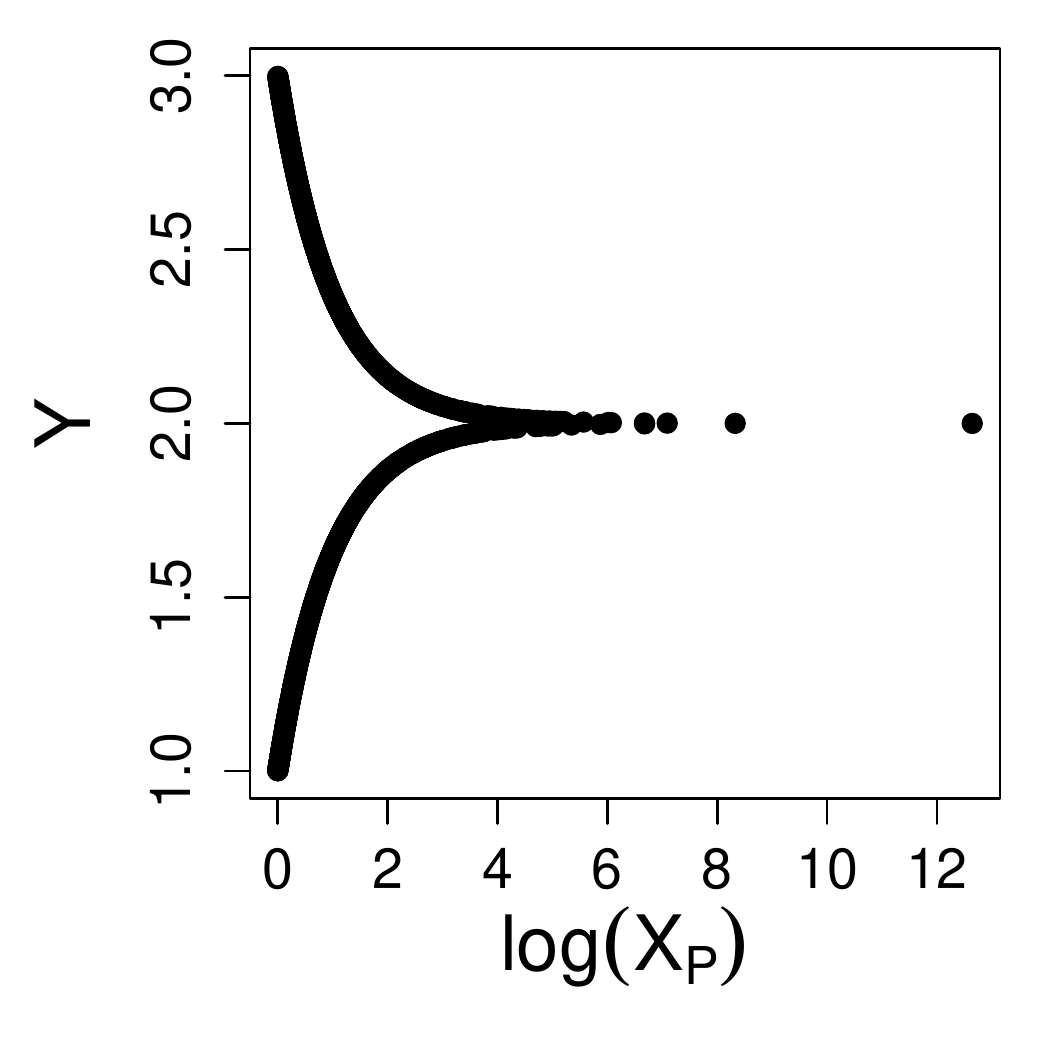}
\hspace{0.5cm}
\includegraphics[width=0.4\textwidth]{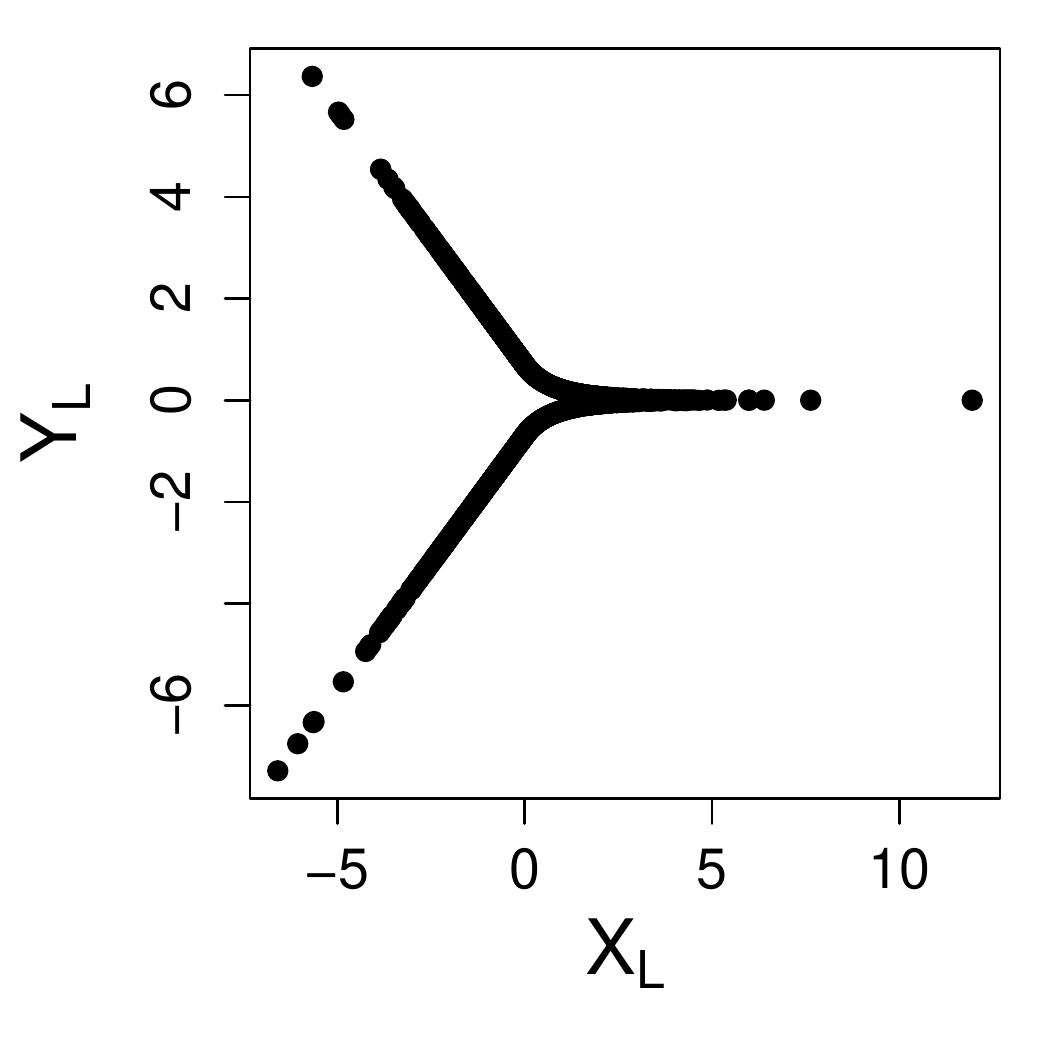}
\caption{Illustration of 2,000 samples for the framework in Example 3.1. The left panel shows the simulated observations $(x_p,y)$ on the original scale, while the right panel corresponds to the transformed samples $(x_L, y_L)$.}
\label{fig:Example31}
\end{figure}

Let $B$ be a discrete random variable that is uniformly distributed on $\{-1,1\}$ and independent of the Pareto distributed random variable~$X_P$. The variable $Y$ is defined as $Y = 2-B/X_P$. For large $X_P$, the values of $Y$ are concentrated around 2 (see Figure~\ref{fig:Example31} left panel). The marginal distribution of $Y$ is $Y\sim\mbox{Uniform}(1,3)$. \cite{Drees2017} present this and the following Example~3.2, to illustrate that the result by \cite{Das2011} linked to the standardisation \eqref{eq:Standardization} of $Y$ does not hold in general.

We again start by transforming the random variables $X_P$ and $Y$ to Laplace margins. Substitution of the values $y$ by their transformed values $y_L$ gives
\[
y = 
\begin{cases}
1 + \exp(y_L) &\mbox{if}~ y_L \leq 0,\\
3 - \exp(-y_L) &\mbox{if}~ y_L > 0.\\
\end{cases}
\]
and the transformation of $X_P$ to Laplace margins is given in \eqref{eq:LaplaceX}. For the case $B=1$, $Y$ takes values smaller than $2$, while only values greater than $2$ are observed for $Y$ when $B=-1$. Therefore, we have to consider the transformation with $y_L \leq 0$ for $B=1$, and $y_L > 0$ for $B=-1$. For the case $B=1$, we find that
\[
y=2-\frac{1}{x} \quad\Leftrightarrow\quad 1+\exp(y_L)=~2-\frac{1}{2}\exp(-x_L)\quad \Leftrightarrow\quad y_L =~\log\left\{ 1 - \frac{\exp(-x_L)}{2}\right\}. 
\]
The final equation implies that we can approximate $y_L$ by $-\frac{1}{2}\exp(-x_L)$ as $x_L\to\infty$. Similar calculations for the case $B=-1$ give that as $x_L\to\infty$
\[
y = 2 + \frac{1}{x}\quad
\Leftrightarrow\quad y_L =~-\log\left\{ 1 - \frac{\exp(-x_L)}{2}\right\} ~\sim~ \frac{1}{2}\exp(-x_L).
\]
Without norming, $Y_L\mid (X_L>x)\rightarrow^P 0$ as $x\rightarrow \infty$. To avoid this degeneracy, we need to take the functions in \eqref{eq:HR20071} to be $a_{\mid X}(x) = 0$ and $b_{\mid X}(x) = \exp(-x)$, the limiting distribution $G_{\mid X}$ then assigns probability $1/2$ to each of the values $z=-0.5$ and $z=0.5$. 

The expression for $b_{\mid X}(x)$ is not of the form for Laplace margins found by \citet{Keef2013a}, with $b_{\mid X}(x)$ tending to zero very rapidly. This form is needed given the speed of convergence of $Y_L\mid(X_L>x)$ towards zero as $x\to\infty$, as seen in Figure~\ref{fig:Example31} right panel. This is not too surprising as it is known that the simple parametric forms of \citet{Keef2013a} for the norming functions do not always hold, with \citet{Papastathopoulos2016} already identifying that
it is possible to have $a_{\mid X}(x) = x\mathcal{L}_a(x)$ and $b_{\mid X}(x) = x^\beta\mathcal{L}_b(x)$, with $\mathcal{L}_a(x)$ and $\mathcal{L}_b(x)$ being slowly varying functions and $\beta\in (-\infty,1)$. Here we have an example that is outside that class with $\beta=0$ and $-\log\{\mathcal{L}_b(x)\}$ being regularly varying. With our norming, the limiting measure $\mu_{Y_L \mid X_L>}$, as defined in \eqref{eq:HR2007}, is 
\[
\mu_{Y_L\mid X_L>}\left( (x_L,\infty] \times [-\infty,y_L] \right) =
\frac{1}{2}\left(\mathbf{I}\left\{-0.5\leq y_L\right\} + \mathbf{I}\left\{0.5\leq y_L\right\}\right) \times \frac{1}{2} \exp(-x_L). 
\]
So, a combination of the standardisation of marginals and random norming, by $X_L$ not $x_L$,
gives a simpler product limit measure than found by \cite{Drees2017}.

\subsection{Example 3.2}
\begin{figure}
\centering
\includegraphics[width=0.4\textwidth]{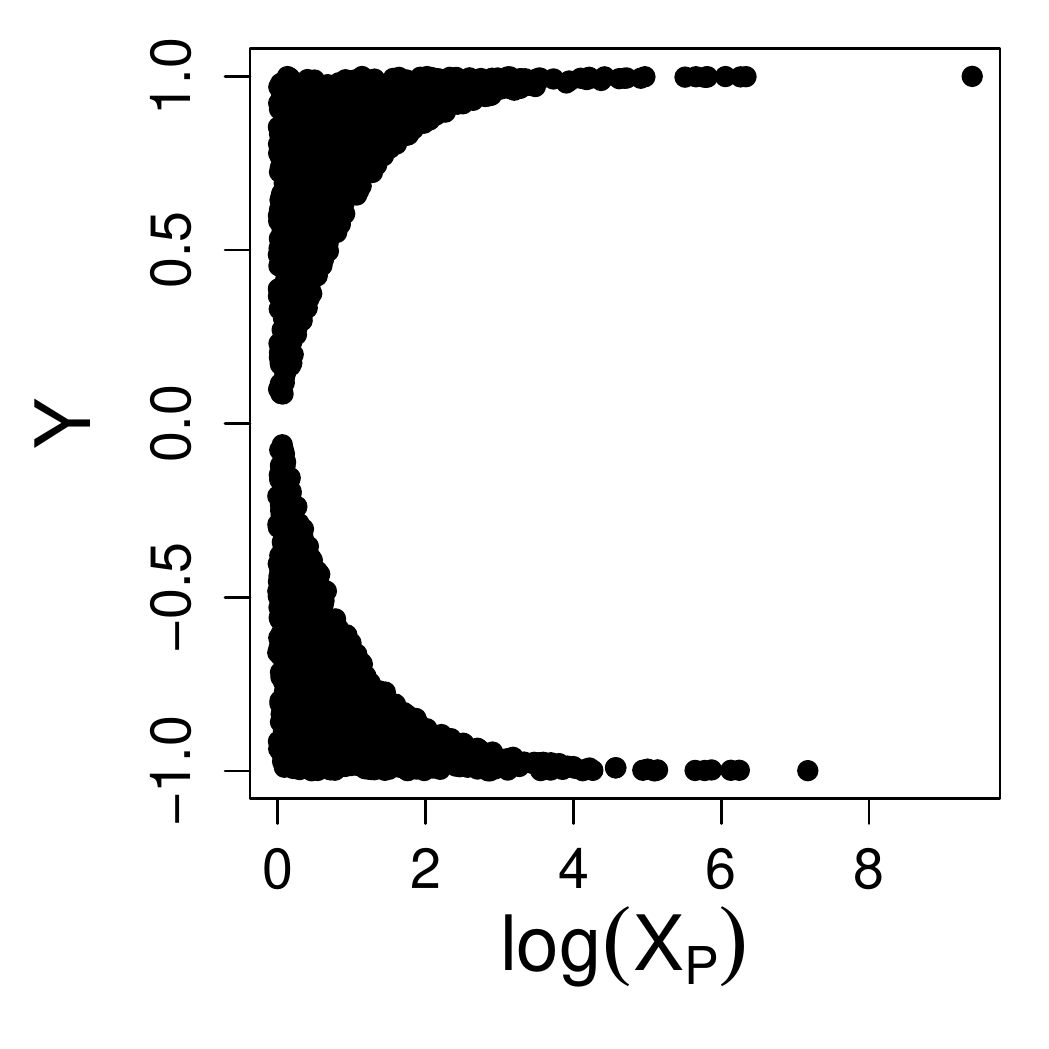}
\hspace{0.5cm}
\includegraphics[width=0.4\textwidth]{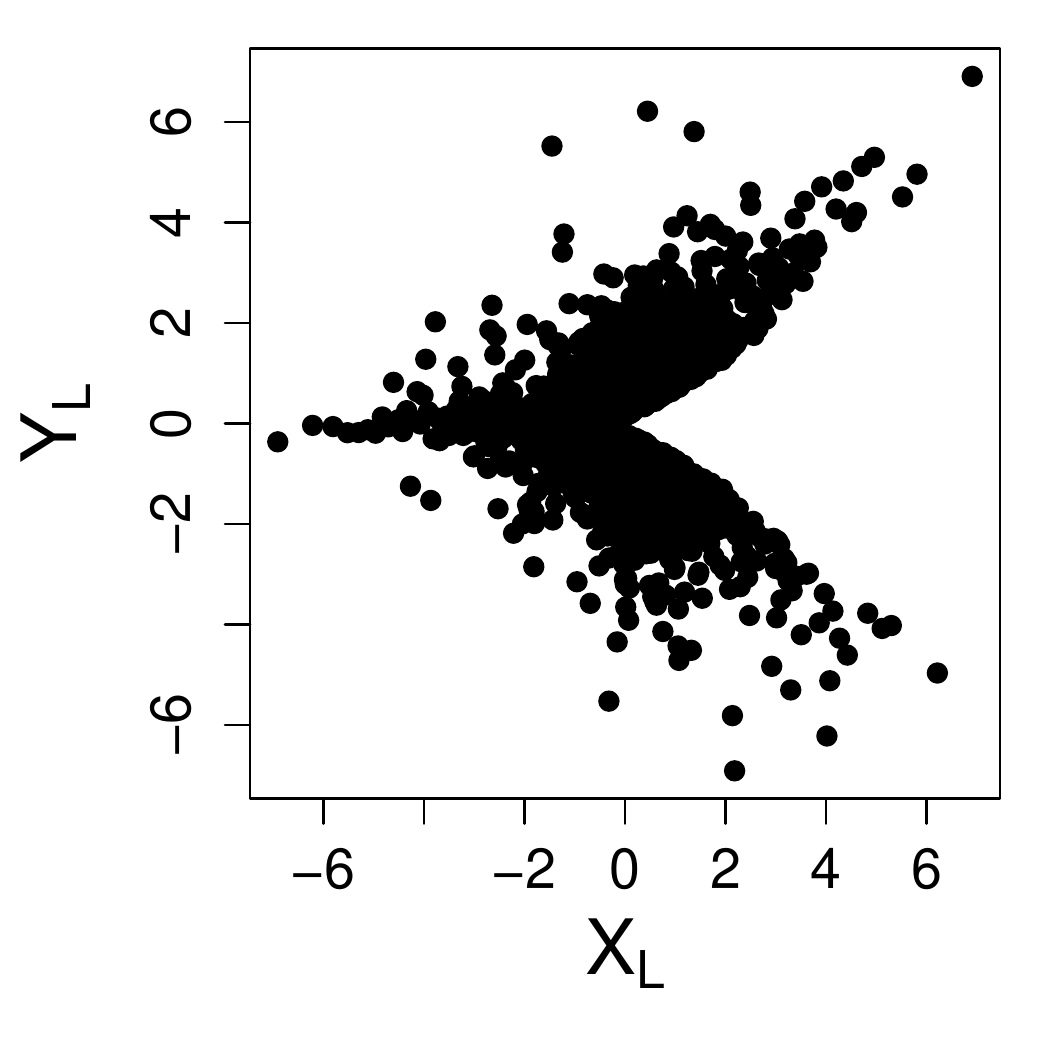}
\caption{Illustration of 2,000 samples for the framework in Example 3.2. The left panel shows the simulated observations $(x_p,y)$ on the original scale, while the right panel corresponds to the transformed samples $(x_L, y_L)$.}
\label{fig:Example32}
\end{figure}

Let $B$ be a discrete random variable that is uniformly distributed on $\{-1,1\}$, $U\sim\mbox{Uniform}(0,1)$, and $X_P$, $B$ and $U$ are all independent. Define $Y = B(1 -U/X_P)$, with the random variable $Y$ taking negative  and positive values for $B=-1$ and $B=1$ respectively. Figure~\ref{fig:Example32} left panel shows that the values of $Y$ are close to $y=-1$ and $y=1$ for large values of $X_P$. For $-1<y<0$, we calculate the marginal distribution of $Y$ as $\Pr(Y<y)= (y+1)\{1-\log(y+1)\}/2$; see Section \ref{sec:Marginal32} for details. Using similar calculations, we find $\Pr(Y<y) = (1+y)/2 + (1-y)\log(1-y)/2$ for $0\leq y<1$.

To transform $Y$ to Laplace margins for $y<0$, which corresponds to $y_L<0$, we use the relationship
$(1/2)(y+1)\left[1-\log(y+1)\right] = (1/2)\exp(y_L)$. Since we cannot find an analytical closed form for $y$ in terms of $y_L$, we consider approximations in order derive the link between $y_L$ and $y$ in the limit as $y\to -1$. The calculations in Appendix~\ref{sec:Limit32} give that
\[
y+1 \sim -\frac{\exp (y_L)}{y_L}
\]
for $y\downarrow -1$. Using similar approximations, we find $1-y \sim  \exp(-y_L)/y_L$ for $y\uparrow 1$. For $B=1$, the limiting behaviour of $Y$, as $X_P$ becomes large, is thus described by
\[
y = - 1 + \frac{u}{x} \quad
\Leftrightarrow\quad -\dfrac{\exp(y_L)}{y_L} - 1 = -1 + \frac{u}{2}\exp(-x_L) \quad
\Leftrightarrow\quad y_L - \log(-y_L) = \log\left(\frac{u}{2}\right) - x_L,
\]
where $u$ denotes the realisation of the random variable $U$. 
Considering $x_L \to \infty$, we obtain that for $B=-1$
\[
y_L = -x_L + \log x_L + o_P(\log x_L), 
\]
where the stochasticity is due to $U$. So, $y_L\overset{p}{\to} -\infty$ as $x_L\to\infty$; this can also be seen in Figure~\ref{fig:Example32} right panel. Similar calculations give $-y_L-\log(y_L)=\log\left(u/2\right)-x_L$ when $B=1$. Hence, for $B=1$, $y_L=x_L-\log x_L+o_P( \log x_L)$ with $y_L\overset{p}{\to}\infty$ as $x_L\to\infty$. 

At first sight these results appear to correspond to there being non-unique choices for $a_{\mid X}$ and $b_{\mid X}$ in \eqref{eq:HT2004} that yield a non-degenerate limiting distribution $G_{\mid X}$. However, there is only one such choice (up to type) with $a_{\mid X}(x) = x$ and $b_{\mid X}(x) = \log x$ ($x>1$) giving $G_{\mid X}$ placing mass of $1/2$ at $\{-\infty\}$ and $\{-1\}$, i.e., $G_{\mid X}(x) = 0.5$ for $-\infty <z<-1$ and $G_{\mid X}(z)=1$ for $-1\le z <\infty$. As in Example~3.1, the derived norming function $b_{\mid X}(x)$ is not of the simple power parametric form of \cite{Keef2013a}. Another possible norming has $b_{\mid X}(x)=x$ with $a_{\mid X}(x)=o(b_{\mid X}(x))$ as $x\rightarrow \infty$, giving $G_{\mid X}$ with mass at of $1/2$ at $\{-1\}$ and $\{1\}$, but this type of norming is not permitted as $b_{\mid X}(x)$ cannot grow as fast as $x$ \citep{Keef2013a}.

\subsection{Example 4.2}
\begin{figure}
\centering
\includegraphics[width=0.4\textwidth]{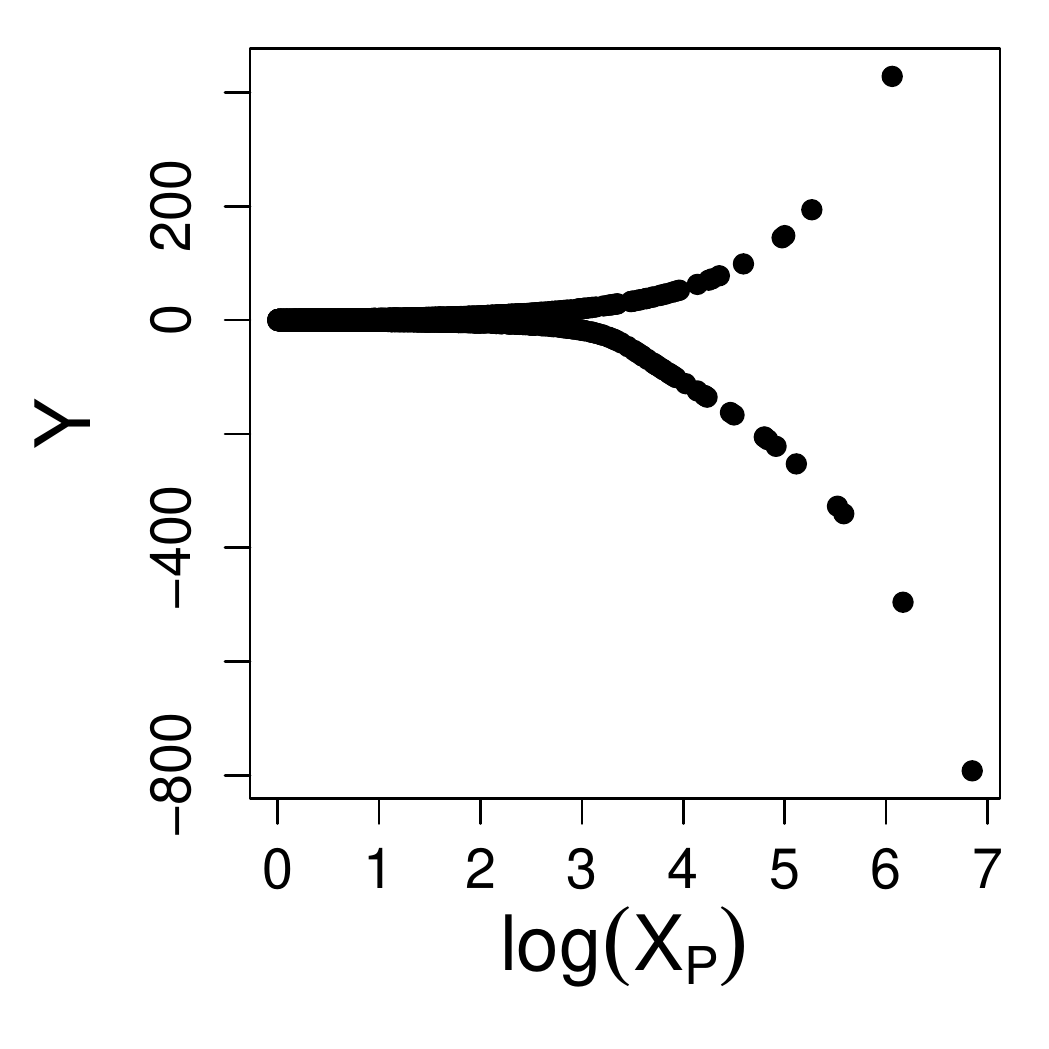}
\hspace{0.5cm}
\includegraphics[width=0.4\textwidth]{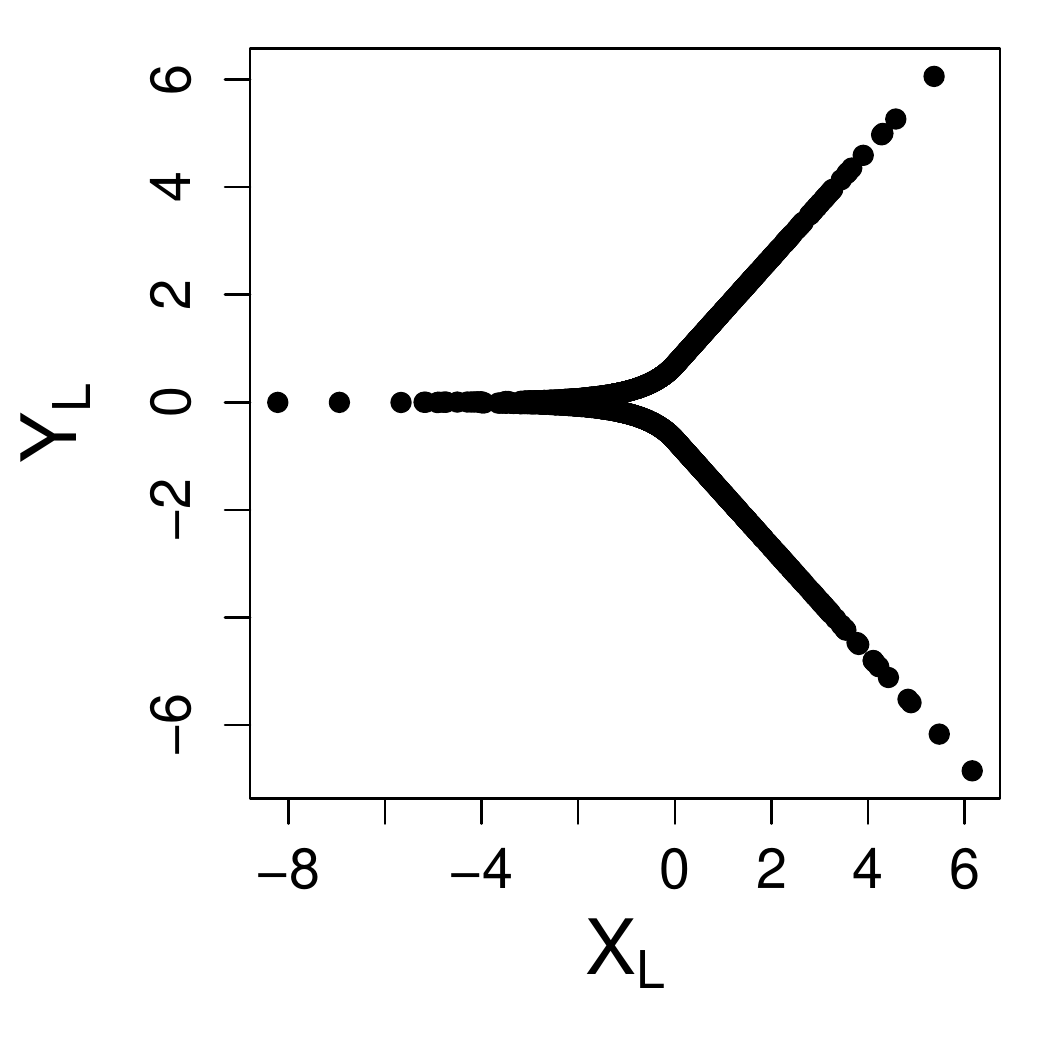}
\caption{Illustration of 2,000 samples for the framework in Example 4.2. The left panel shows the simulated observations $(x_p,y)$ on the original scale, while the right panel corresponds to the transformed samples $(x_L, y_L)$.}
\label{fig:Example42}
\end{figure}

Let $B$ be a discrete random variable that is uniformly distributed on $\{0,1\}$. Define the function $g(x):=x(2 + \sin\log x)$ for $x\geq 1$ and we consider $Y = B X_P + (1-B)\left\{-g^{-1}\left(2 X_P\right)\right\}$, where $g^{-1}$ is the inverse of $g$. Figure~\ref{fig:Example42} indicates that $Y$ tends to $-\infty$ and $+\infty$ as $X_P$ becomes large. The purpose of this example in \cite{Drees2017} is to illustrate that ($X_P,Y$) being multivariate extreme value distributed is not a sufficient condition for $Y\mid (X_P>t)$, as $t\to\infty$, to lie in the class of CEVMs of the form \eqref{eq:HR2007}. 

We derive the marginal distribution of $Y$ as 
\[
\Pr(Y<y) = 
\begin{cases}
1 / g(-y)  &\mbox{if}~y<-1,\\
1/2      &\mbox{if}~-1\le y\le 1\\
1 - 1 / 2y & \mbox{if}~y>1,
\end{cases}
\]
and the calculations are provided in Appendix A.3. Transformation of $Y$ to Laplace margins gives $y = -g^{-1}\left\{2/\exp(y_L) \right\}$ for $y\leq -1$, and $y = \exp(y_L)$ when $y\geq 1$. The limiting behaviour of the transformed variable $Y_L$ as $X_L$ becomes large, and for $B=1$, is given by $\exp (y_L) \sim 2 \exp (x_L)$, which is equivalent to $y_L = \log 2 + x_L+o(1)$ as $x_L\rightarrow \infty$. For the case $B=0$, we find $y_L = -\log 2 - x_L+o(1)$ as $x_L\rightarrow \infty$.  This symmetry in the limiting behaviour on Laplace marginal distributions is also visible in Figure~\ref{fig:Example42} right panel. Consequently, we have that as $x_L\rightarrow \infty$,
\[
Y_L =
\begin{cases}
\log 2 + X_L +o(1) & \mbox{if}~B=1,\\
-\log 2 - X_L +o(1) & \mbox{if}~B=0.\\
\end{cases}
\]

Defining $a_{\mid X}(x) = x$ and $b_{\mid X}(x)=1$ yields a non-degenerate limiting distribution $G_{\mid X}$ with $G_{\mid X}(z) = 0.5$ for $-\infty < z<\log 2$ and $G_{\mid X}(z) = 1$ for $\log 2\le z < \infty$. While the normalising functions $a_{\mid X}(x) = -x$ and $b_{\mid X}(x)=1$ also yield a non-degenerate limiting distribution $G_{\mid X}$, this choice is not permissible because $G_{\mid X}$ would have mass at $+\infty$ \citep{Keef2013b}. Consequently, the normalising functions $a_{\mid X}$ and $b_{\mid X}$ are well-defined (up to type). Furthermore, this example shows that a transformation to Laplace margins can result in the distribution of $Y_L\mid (X_L>t)$ as $t\to\infty$ being in the class of conditional extreme models by \cite{Keef2013a}, although $Y\mid (X_P>t)$ as $t\to\infty$ does not lie in the class of CEVMs introduced by \cite{HeffernanResnick2007}. When we consider the distributions of $X_L \mid (Y_L=y_L)$, where $y_L>0$, there is deterministic relationship between $X_L$ and $Y_L$, and thus there cannot exist a non-degenerate limiting distribution $G_{\mid Y}$, but the behaviour is trivial $X_L=Y_L-\log 2 \mid Y_L>y_L$ for any $y_L>0$. 

\subsection{Example 4.4}

Define the function $g_c(u) = u(1 + c \sin\log u)$, where $0<u\le 1$ and $|c|<1/\sqrt{2}$, and $\psi_c(z) = g_c^{-1}(1/z)$ with $z\geq 1$. Let $Z_P$ be standard Pareto distributed, $\Pr(Z_P>z)=1-z^{-1}~(z>1)$, and $B$ be discrete and uniformly distributed on $\{1,2,3,4\}$ and also independent of $Z_P$. The random variables $X$ and $Y$ are then defined as
\begin{equation}
(X,Y) := 
\begin{cases}
(2- \psi_{1/2}(Z_P),~2-1/Z_P )&\mbox{if}~B=1,\\
(2- \psi_{-1/2}(Z_P),~2-1/\sqrt{Z_P} )&\mbox{if}~B=2,\\
(1- 1/Z_P,~2-1/Z_P )&\mbox{if}~B=3,\\
(2- 1/Z_P,~1-1/Z_P )&\mbox{if}~B=4.
\end{cases}
\label{eqn:Ex44setup}
\end{equation}
The purpose of this example by \citet{Drees2017} is to show that ($X,Y$) does not lie in the class of multivariate extreme value models despite $Y\mid (X>t)$ and $X\mid (Y>t)$ belonging to the CEVM class by \cite{HeffernanResnick2007} as $t\to2$. This inconsistency of the CEVM class with $(X,Y)$ being in the domain of attractions of a  bivariate extreme value distribution indicates that these conditional distributions fall outside the framework of the standard assumptions for bivariate extreme values.

\begin{figure}
\centering
\includegraphics[width=0.4\textwidth]{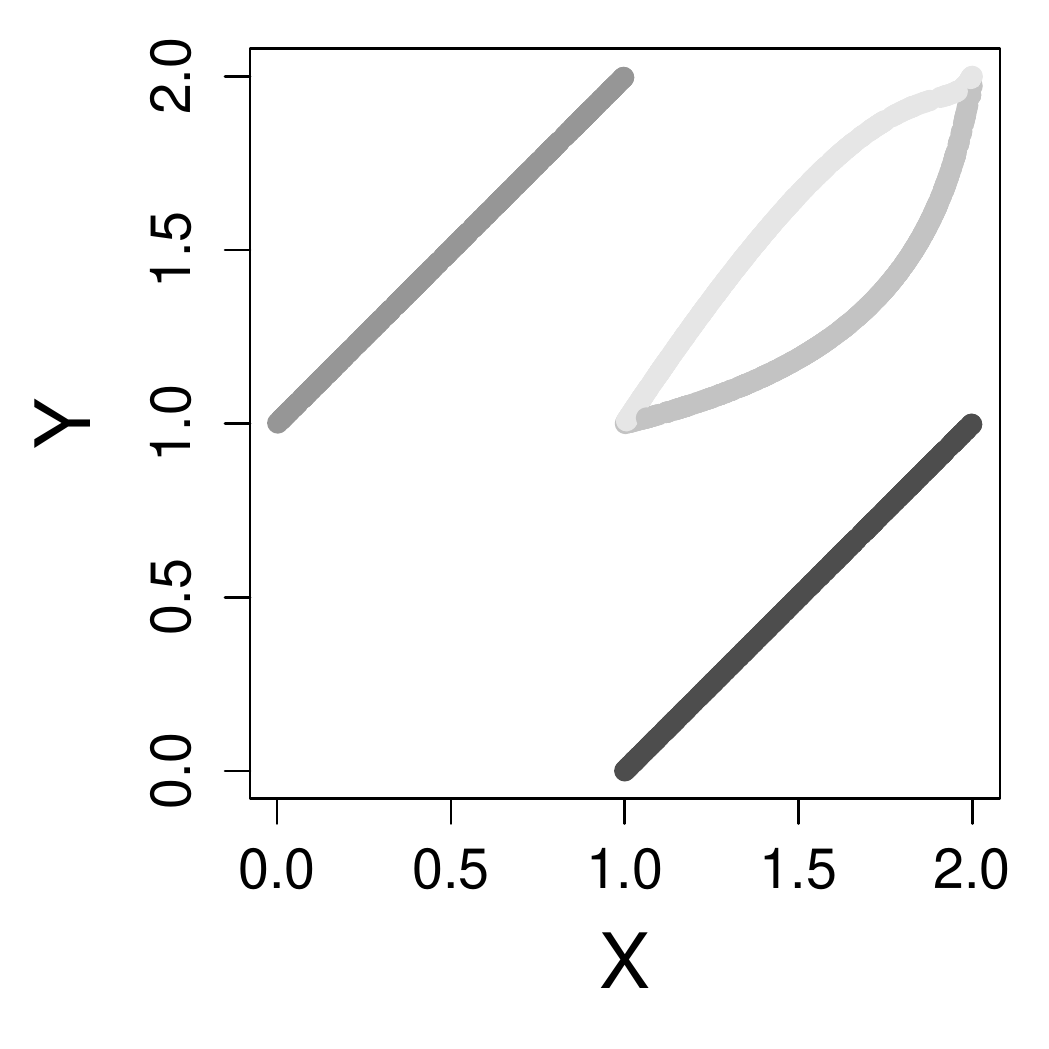}
\hspace{0.5cm}
\includegraphics[width=0.4\textwidth]{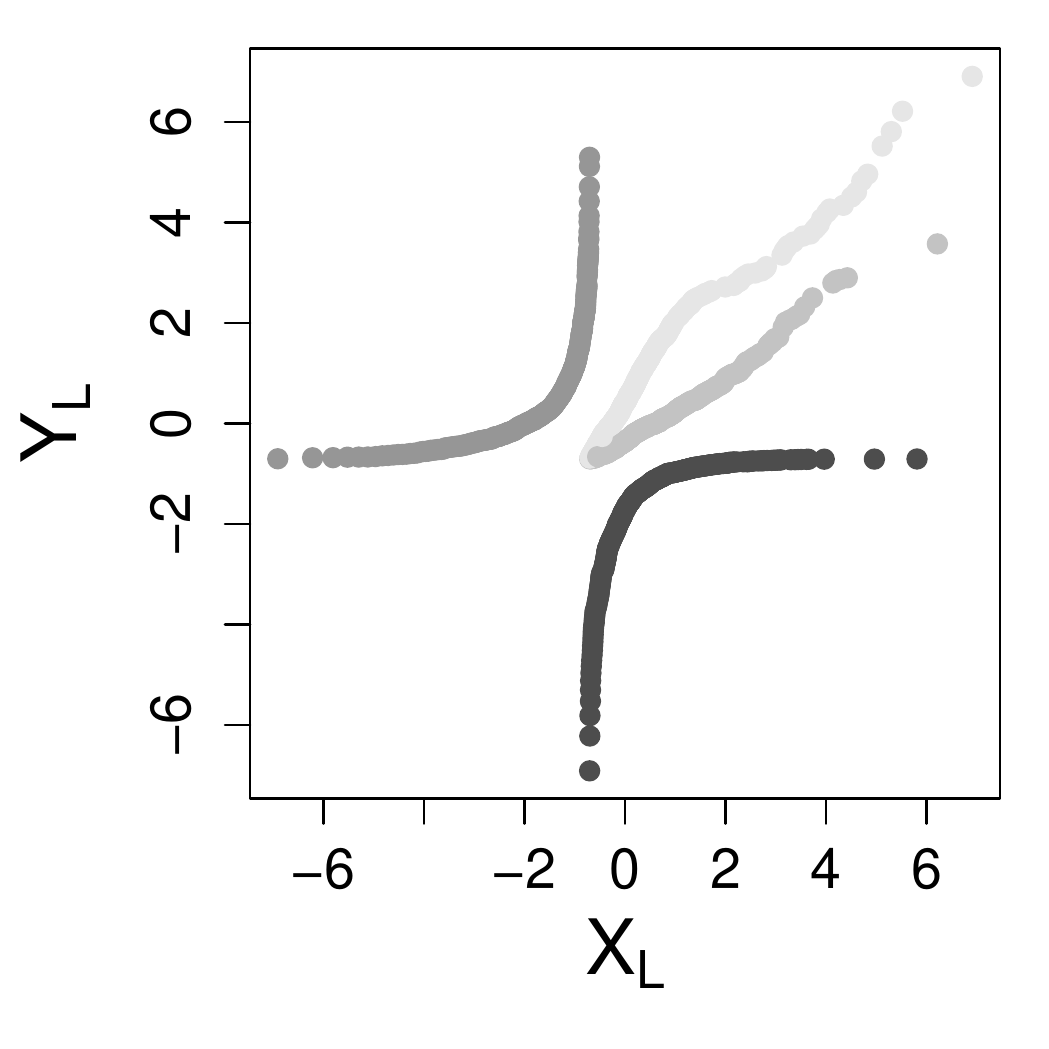}
\caption{Illustration of 2,000 samples for the framework in Example 4.4. The left panel shows the simulated observations $(x,y)$ on the original scale, while the right panel corresponds to the transformed samples $(x_L, y_L)$, with the behaviour for each value of $B$ highlighted: the lightest shade corresponds to $B=1$, and the points with the darkest shade are the samples for $B=4$.}
\label{fig:Example44}
\end{figure}

We now investigate this inconsistency for the bivariate distribution~\eqref{eqn:Ex44setup} after marginal standardisation. We start by calculating the marginal distributions of $X$ and $Y$ to see if they are individually in the domain of attraction of the univariate extreme value distribution~\eqref{eq:Tail}. Figure \ref{fig:Example44} left panel shows that the random variable $X$ ($Y$) respectively can take values between 0 and 1 when $B=3$ ($B=4$), while $B\neq3$ ($B\neq 4$) leads to the values of $X$ ($Y$) lying between 1 and 2. The cumulative distribution function of $X$ is
\[
\Pr(X\leq x ) = 
\begin{cases}
x/4 & \mbox{if}~0\leq x\leq 1,\\
3x/4 - 1/2 &\mbox{if}~1\leq x\leq 2,
\end{cases}
\]
and for $Y$ we have
\[
\Pr(Y\leq y ) = 
\begin{cases}
y/4 & \mbox{if}~0\leq y\leq 1,\\
y/2 - (2-y)^2/4 &\mbox{if}~1\leq y\leq 2.
\end{cases}
\]
Detailed calculations for $\Pr(X\leq x )$ and $\Pr(Y\leq y )$ are provided in Appendix A.4. For these two marginals it is straightforward to show that they are each in the domain of attraction of the univariate extreme value distribution with parameters $\gamma_X=\gamma_Y=-1$.

Transformation of $X$ and $Y$ to Laplace margins gives
\[
x = 
\begin{cases}
2 \exp(x_L)         &\mbox{if}~x_L\leq-\log 2,\\
2/3 + (2/3) \exp(x_L) &\mbox{if}~-\log 2<x_L\leq0,\\
2 - (2/3) \exp(-x_L)  &\mbox{if}~x_L>0,
\end{cases}
\]
and 
\[
y = 
\begin{cases}
2 \exp(y_L)              &\mbox{if}~y_L\leq-\log 2,\\
3 - \sqrt{5-2\exp(y_L)}  &\mbox{if}~-\log 2<y_L\leq0,\\
3 - \sqrt{1+2\exp(-y_L)} &\mbox{if}~y_L>0.
\end{cases}
\]
Applying these transformations does not change the property that the marginal distributions are in the domain of attractions of univariate extremes value distribution, only now $\gamma_{X_L}=\gamma_{Y_L}=0$, and we still have that $(X_L,Y_L)$ is not in domain of attraction of a bivariate extreme value distribution. The following calculations show that, with standardisation to Laplace margins, the conditional limiting distribution of $Y_L\mid(X_L>x_L)$ fails to meets the conditions of \citet{HT2004} as $x_L\to\infty$, unlike the conditional $Y\mid(X>t)$, as $t\to2$, that falls in the CEVM class of \citet{HeffernanResnick2007}.

We explore the conditional distributions by looking at the relations between $(X_L,Y_L)$ for each value of $B$. For $B=1$, the expressions $X = 2-\psi_{1/2}(Z)$ and $Y=2-1/Z$ give $Y = 2 - g_{1/2}(2-X)$. To study the limiting behaviour, we again replace $X$ and $Y$ by their Laplace distributed transformed expressions. Figure~\ref{fig:Example44} right panel shows that, for $B=1$, large values of $X_L$ lead to large values of $Y_L$ and we thus consider the equality
\begin{align*}
3 - \sqrt{1 + 2\exp(-y_L)} &= 2 - g_{1/2} \{(2/3) \exp(-x_L)\}\\
&= 2 - (1/3) \exp(-x_L)[2 + \sin\{ \log (2/3) - x_L \}].
\end{align*}
For notational brevity, we define $h_1(x_L) = (1/3) \left[2 + \sin\{\log (2/3) - x_L\}\right]$ and we note $(1/3) \leq h_1(x_L) \leq 1$. By simplifying the terms and taking squares on both sides, we get
\[
1+2\exp(-y_L) = 1 + 2\exp(-x_L) h_1(x_L) + \exp(-2x_L) \left[h_1(x_L)\right]^2.
\]
Further simplifying the terms and taking logs, we end up with
\begin{equation}
y_L = x_L - \log h_1(x_L) + \log\left[ 1 + (1/2) \exp(-x_L) h_1(x_L)\right],
\label{eqn:diagonal.asypmt}
\end{equation}
which we can write as $y_L = x_L - \log h_1(x_L) + (1/2) \exp(-x_L) h_1(x_L) + O\left(\exp(-2x_L)\right)$.

For $B=2$, the expressions $X=2-\psi_{-1/2}(Z_P)$ and $Y=2-1/Z_P$ give that $Y = 2 - \sqrt{g_{-1/2}(2-X)}$. On Laplace scale, we then have
\[
3 - \sqrt{1 + 2\exp(-y_L)} 
~=~ 2 - \sqrt{g_{-1/2} \left\{(2/3) \exp(-x_L)\right\}}\\
~=~ 2 - \exp(-x_L/2) h_2(x_L),
\]
where $h_2(x_L) = \sqrt{ (2/3)  - (1/3) \sin\{ \log (2/3) - x_L\}}$. By taking squares on both sides,
\[
1+2\exp(-y_L) = 1+2\exp(-x_L/2) h_2(x_L)+ \exp(-x_L)\left[h_2(x_L)\right]^2.
\]
Following the same steps as for $B=1$ yields
\[
y_L = x_L/2 - \log h_2(x_L) + (1/2)\exp(-x_L/2) h_2(x_L) + O\left(\exp(-x_L)\right).
\]

The case $B=3$ leads to $x_L < -\log 2$, i.e., $x_L$ is not becoming large, and thus this mixture component can be ignored when studying $Y_L \mid (X_L>x_L)$ as $x_L\to\infty$. Finally, $B=4$ gives $Y = X-1$ with
$2\exp(Y_L) = 2 - (2/3)\exp(-X_L) - 1$. Consequently, 
\[
y_L = \log\left[1/2 - (1/3)\exp(-x_L)\right]~\sim~ -\log(2)\qquad\mbox{ as}~x_L\to\infty.
\]

Now we consider combining the different mixture components and we set $a_{\mid X}(x) = x - \exp(-x) h_1(x)$ and $b_{\mid X}(x) = -\log h_1(x) + (1/2) h_1(x) \exp(-x)$ in \eqref{eq:HR20071}; as $-\log h_1(x)\geq 0$ and $h_1(x)> 0$ this gives $b_{\mid X}(x)>0$ as required. The limiting behaviour of $Y_L\mid(X_L>x_L)$ as $x_L\rightarrow \infty$, for $B=1$, is then
\[
\frac{Y_L - X_L +\exp(-X_L) h_1(X_L)}{-\log h_1(X_L) + (1/2) h_1(X_L) \exp(-X_L)} \sim 1. 
\]
For the remaining components, $B=2$ and $B=4$, $\lim_{x_L\to\infty}(Y_L-a_{\mid X}(X_L))/{b_{\mid X}(X_L)}=-\infty$ for $\{X_L>x_L\}$. However, there is no limiting distribution $G_{\mid X}$ as in \eqref{eq:HR20071} because $\Pr(B=1\mid X_L>x_L)$ oscillates between $1/6$ and $1/2$ as $x_L\to\infty$, that is, $\Pr\left( \{Y_L-a_{\mid X}(x_L)\} / b_{\mid X}(x_L) \le z \mid X_L > x_L \right)$ does not converge. This oscillating behaviour is found by considering $\Pr(B=1\mid X>t)$ for $1<t<2$. Using similar calculations as in Appendix A.4, we find
\[
\Pr(B=1\mid X>t) = \frac{\Pr(X>t\mid B=1) \Pr(B=1)}{\Pr(X>t)} = \frac{1}{3}+\frac{1}{6}\sin\log(2-t),
\]
which oscillates between 1/6 and 1/2 as $t \to 2$, and this implies that $\Pr(B=1\mid X_L>x_L)$ oscillates between 1/6 and 1/2 as $x_L\to\infty$. Consequently, $Y_L\mid(X_L>x_L)$ as $x_L\to\infty$ does not fall in the class of conditional extreme value models by \cite{Keef2013a}, despite $Y\mid(X>t)$, as $t\to 2$, being in the CEVM class by \citet{HeffernanResnick2007}, see \cite{Drees2017}.

So far we have focused on the conditional distribution of $Y_L\mid (X_L>x_L)$ for $x_L\rightarrow \infty$, but there is  also interest in the asymptotic behaviour of the reverse conditional $X_L\mid (Y_L>y_L)$ as $y_L\to\infty$. Here Figure~\ref{fig:Example44} provides some insight into what happens, with only the mixture terms corresponding to $B=1$ and $B=3$ contributing to the tail of $Y_L$. Further, expression~\eqref{eqn:Ex44setup} shows that $Y_L\mid (B=1)$ and $Y_L \mid (B=3)$ are identical, which gives that the limiting distribution of $X_L \mid (Y_L>y_L)$ as $y_L\rightarrow \infty$ must be a mixture distribution with weights $1/2$ on each component. When $B=3$ we see that $X_L$ does not grow with $Y_L$, so with any norming on $X_L$ that is required to handle the growth of $X_L$ with $y_L$ in $X_L\mid (Y_L>y_L)$ will lead to mass tending to $-\infty$ when $B=3$. So it remains to consider the $B=1$ case. The deterministic relationship in expression~\eqref{eqn:diagonal.asypmt} between $X_L$ and $Y_L$ gives $Y_L>X_L$ conditional on $X_L$ being above a sufficiently high threshold, because $-\log h_1(x) \geq 0$ and $\log[1+(1/2) \exp(-x) h_1(x)]>0$ for all $x$. Furthermore, the relation between $X_L$ and $Y_L$ is bijective, because the first derivative in~\eqref{eqn:diagonal.asypmt} is strictly positive. Consequently, we can invert the relation between $Y_L$ and $X_L$ when $B=1$ and this gives for $y_L\to\infty$ that 
\[
x_L= y_L-Q(y_L),
\]
where $Q(y_L)>0$ is an oscillator function that is bounded above. When $B=1$, we thus obtain the limiting behaviour
\[
\frac{x_L - y_L}{Q(y_L)} \sim -1. 
\]
Hence we have that for all $z\in \mathbb{R}$, as $y_L\rightarrow \infty$
\[
\Pr\left(\left.\frac{X_L - Y_L}{Q(Y_L)}<z\right| Y_L>y_L\right)\rightarrow 0.5[1+\mathbf{I}(z>-1)].
\]
Thus, the reverse conditional has a more straight-forward behaviour.

We further note that the transformation to Laplace margins does not lead to ($X_L,Y_L$) lying in the class of multivariate extreme value models. \cite{Drees2017} showed that ($X,Y$) is not multivariate extreme value distributed either. Consequently, this is an example for which the limiting behaviour $Y_L \mid (X_L>x_L)$ is not in the class of conditional extremes models by \cite{Keef2013a} and $(X_L,Y_L)$ does not lie in the domain of attraction of a multivariate extreme value distribution.

\section{Implications to the conditional extremes model}
\label{sec:Implications}

The conditional extreme modelling frameworks by \cite{HT2004} and \cite{Keef2013a} have been widely adopted in practice. In an effort to introduce these methods to an even wider class of statistical problems, and to establish a stronger link to other multivariate extremes models, \cite{HeffernanResnick2007}, amongst others, considered settings with non-standardized margins. The counterexamples in \cite{Drees2017} demonstrate, however, that the link between this class of conditional extremes models and other multivariate extremes models is quite subtle and not as extensive as described by, for instance, \cite{Das2011} and \cite{Resnick2014}.

In Section~\ref{sec:Examples}, we showed that some of the issues, raised by Examples~2.3 to 4.4 in \cite{Drees2017}, can be resolved by standardisation to common Laplace margins. As such, it is worth discussing the practical implications of our results, with a particular focus regarding the application of the \cite{HeffernanResnick2007} and \cite{Keef2013a} conditional extreme models.

As stated in the Introduction, we believe that working with completely different marginal tail behaviours imposes a major restriction on a conditional extremes approach using affine transformations. Our calculations clearly demonstrate that standarisation to common margins may have some benefits. In Example~2.3, this fixed choice of standardisation implied a unique limit measure $G_{\mid X}$, while the CEVM framework by \cite{HeffernanResnick2007} allowed the limit measure $\mu_{Y\mid X>}$ to vary with the standardisation used. This example also highlighted that it is necessary to allow $G_{\mid X}$ to have mass at $\{-\infty\}$ in the \cite{Keef2013a} framework, because the limit measure might otherwise be degenerate. Given these findings, we believe our results illustrate the versatility of the \citet{HT2004} conditional multivariate extremes framework.

Nevertheless, we note that the examples of \cite{Drees2017} illustrate some statistical limitations of the \cite{HT2004} framework even with standardised Laplace marginals of \citet{Keef2013a}. Two particular areas relate to handling mixture distributions for $G_{\mid X}$ and the choice of parametric families for the normalising functions $a_{\mid X}$ and $b_{\mid X}$. We discuss these in turn below.

Many of the examples of \cite{Drees2017} involved a mixture structure for $(X,Y)$, and hence also for $(X_L,Y_L)$. Although it was possible to identify normalising functions to give a non-degenerate $G_{\mid X}$ in these cases, it was no surprise that $G_{\mid X}$ was also a mixture distribution. From a statistical perspective the only complication with $G_{\mid X}$ having a mixture structure is when $G_{\mid X}$ puts an atom of mass at $\{-\infty\}$; with Example~3.1 being the only example where $\{-\infty\}$ has mass zero. The complication with limiting mass at $\{-\infty\}$  is that at non-asymptotic levels of $x_L$ this mass will be at a finite value with its precise value depending on the associated conditioning value, e.g., $x_L$ in this set up. Statistical methods have recently been developed by \cite{Tendijck2021} which extend the Heffernan–Tawn conditional extreme value model for handling exactly this situation.

\citet{Keef2013a} propose parsimonious canonical parametric families for the normalising functions $a_{\mid X}(x)=\alpha x$ and $b_{\mid X}(x)= x^\beta$ which appears suitable for a wide range of published  data applications. The examples in \cite{Drees2017} add to the list (first noted by \citet{Papastathopoulos2016}) of theoretical joint distributions for $(X_L,Y_L)$ with normalising functions that lie outside the canonical class. Clearly, the canonical families cannot be extended to cover all of these theoretical examples in a parsimonious way. So the most natural line of future research is to identify if it is possible to quantify the errors that can  arise from the inappropriate usage of the canonical families in these examples, with the error relating to the bias of estimated  probabilities of extremes events for finite extrapolations.

\section{Discussion}
\label{sec:Discussion}

By investigating Examples~2.3 to 4.4 by \cite{Drees2017} throughout this paper, we found some interesting differences between the conditional extremes frameworks by \cite{Keef2013a} and \cite{HeffernanResnick2007}. In Example~2.3, the standardisation to Laplace margins implied a unique limit measure $G_{\mid X}$, while the CEVM framework by \cite{HeffernanResnick2007} allowed the limit measure $\mu_{Y\mid X>}$ to vary with the standardisation used. This example also highlighted that it is necessary to allow $G_{\mid X}$ to have mass at $\{-\infty\}$, because the limit measure might otherwise be degenerate. Consequently, while \cite{Drees2017} advocate for the condition $\mu_{Y\mid X>}\left(\{-\infty,\infty\}\times E^{(\gamma_X)}\right)=0$ to ensure uniqueness, in our calculations it was sufficient to require $\lim_{z\to\infty} G_{X}(z)=1$ for the conditional extremes model by \cite{Keef2013a}.

A non-degenerate limit measure $G_{\mid X}$ was found in Examples 3.1 and 3.2, however, the functions $a_{\mid X}$ and $b_{\mid X}$ in \eqref{eq:HR2007} were not of the simple parametric form of \cite{Keef2013a}. As mentioned in Section~\ref{sec:Implications}, it was already known that the canonical families by \cite{Keef2013a} cannot cover all cases in a parsimonious way, and our results add to this set of examples. Example 3.1 further shows that standardisation to Laplace margins can result in a non-degenerate limit measure, despite there not existing a standardisation of the form by \cite{Das2011} for the CEVM framework. Conversely, Example~4.4 shows that a model may belong the CEVM class by \cite{HeffernanResnick2007}, but not fall in the class of conditional extreme models by \cite{Keef2013a}.

Finally, Example~4.2 in \cite{Drees2017} showed that $(X,Y)$ being multivariate extreme value distributed is not sufficient for the distributions of $Y\mid (X>t)$ and $X \mid (Y>t)$ to be in the CEVM class as $t$ approaches the upper end of $X$ and $Y$ respectively, while Example~4.4 illustrated that $X\mid(Y>t)$ and $Y\mid (X>t)$ being CEVM does not imply that the distribution of $(X,Y)$ is in the domain of attraction of multivariate extreme value distributions. However, after standardisation of $(X,Y)$ to Laplace margins $(X_L,Y_L)$, our working shows that Examples 4.2 and 4.4 are ruled out as being evidence for the limit~\eqref{eq:HR20071}, and the associated result for $X_L\mid Y_L$, not being equivalent to the domain of attraction condition of a bivariate extreme value distribution. Consequently, the link between the conditional extremes models by \cite{Keef2013a} and the class of multivariate extreme value distributions remains an open research question. 

%\CR{Jon: Please add something to address the following referee comment:\\
%"I would also like to point out that there is a fundamental difference in the models studied in this article and the CEV as studied in \cite{Drees2017}. In the latter article, a non-random norming was used as in Equation (6) of the present article (the functions $f$ and $g$ to standardize $Y$ are applied to the non-random $t$, as $t\to\infty$). Yet, in Keef et al. (2013a) and in the analysis of the examples in the present article a random norming as in Equation 5 of the present article is used. I could not find a proper discussion of this difference in the text but it can be essential. In particular, a random norming would usually lead to a product limit (cf. Example 3.1 which has a product limit in contrast to the original example). So, while I haven't looked into the details, I assume that the counterexamples are, at least to some extent, not overcome by standardisation but by random norming."}

\appendix

\section{Appendix: Technical details}

\subsection{Marginal distribution in Example 3.2}
\label{sec:Marginal32}

For $-1 < y < 0$, we have 
\begin{align*}
\Pr(Y<y) 
& ~=~ \frac{1}{2}\,\Pr\left(-1 + \frac{U}{X_P} < y\right)\\
& ~=~ \frac{1}{2}\int_0^1 \Pr\left(X_P > \frac{u}{y+1}\right) f_U(u) \mathrm{d}u\\  
&~=~ \frac{1}{2}\int_0^{y+1} \Pr\left(X_P > \frac{u}{y+1}\right) \mathrm{d}u + \frac{1}{2}\int_{y+1}^1 \Pr\left(X_P > \frac{u}{y+1}\right) \mathrm{d}u\\  
&~=~ \frac{1}{2}\int_0^{y+1} 1\, \mathrm{d}u + \frac{1}{2}\int_{y+1}^1 \frac{y+1}{u} \mathrm{d}u\\  
&~=~ \frac{1}{2}(y+1)\{1-\log(y+1)\}.
\end{align*}

\subsection{Approximation of the limiting behaviour in Example 3.2}
\label{sec:Limit32}

Our aim is to find an approximation of the relation between $y$ and the transformed value $y_L$ as $y\to -1$. We first note that $(y+1)\left[1-\log(y+1)\right]\approx -(y+1)\log(y+1)$ for $y$ close to $-1$. The next step is to define $y+1 = \epsilon \exp(y_L)$ for some $\epsilon>0$. If we set $\epsilon = -1/y_L$, we obtain 
\[
-\frac{1}{y_L}\exp(y_L)\left[ -\frac{1}{y_L} - \exp(y_L)\right]~=~ \exp(y_L).
\]
Taking logs on both sides,  we get $-\log y_L + y_L + \log \left( y_L - \log y_L \right)~=~ y_L$ and the difference of the two sides becomes negligible for $y_L$ large enough. Consequently, we have for $y\approx-1$ that
\[
y \approx -\frac{\exp (y_L)}{y_L} - 1.
\]

\subsection{Marginal distribution in Example 4.2}
\label{sec:Marginal42}

From the definition of $Y= BX_P + (1-B) \{-g^{-1}(2X_P)\}$, it is clear that $Y$ can only take negative values for $B=0$. Further, $g^{-1}: [2,\infty) \to [1,\infty)$ gives that $-g^{-1}(2X_P)$ will only take values smaller than -1. Consider $y\leq -1$ and we calculate
\begin{align*}
\Pr(Y \le y ) &= \frac{1}{2}\Pr\left(-g^{-1}(2X_P) \le y\right)\\
&=\frac{1}{2} \Pr\left ( X_P > \frac{1}{2} g(-y)\right) \\
&=\frac{1}{g(-y)}.
\end{align*}
Finally, the case $B=1$ results in $Y$ only taking values greater than 1 with
\[
\Pr(Y\leq y) = \frac{1}{2} + \frac{1}{2} P( X_P \le y) = \frac{1}{2} + \frac{1}{2}\left(1- \frac{1}{y}\right) = 1- \frac{1}{2y}
\]
as described in Example 4.2.

\subsection{Marginal distributions in Example 4.4}

Only $B=3$ leads to a value of $X$ between 0 and 1, while $X$ lies between 1 and 2 for $B=1,2,4$. Considering $0<x<1$,
\[
\Pr(X\le x) = \frac{1}{4} \Pr\left(1-\frac{1}{Z_P} \le x\right) = \frac{1}{4} \Pr\left(Z_P \le x-1 \right) = \frac{1}{4}\left[1-(x-1)\right] = \frac{x}{4}.
\]
When $1<x<2$, $\Pr(X\le x)$ is given by
\begin{align*}
&\frac{1}{4}\left[1 + \Pr(X\le x\mid B=2) + \Pr(X\le x\mid B=3) + \Pr(X\le x\mid B=4)\right]\\
&=\frac{1}{4}\left[1 + \Pr\left(2-\psi_{1/2}(Z_P)\le x\right) + \Pr\left(2-\psi_{-1/2}(Z_P)\le x\right) + \Pr\left(2-\frac{1}{Z_P}\le x\right)\right] \\
&= \frac{1}{4}\left[1 + \Pr\left\{2-x\le g^{-1}_{1/2}\left(\frac{1}{Z_P}\right)\right\} + \Pr\left\{2-x\le g^{-1}_{-1/2}\left(\frac{1}{Z_P}\right)\right\} + \Pr\left(Z_P\le \frac{1}{2-x}\right)\right]\\
&= \frac{1}{4}\left[1 + \{1 - g_{1/2}(2-x)\} + \{1 - g_{-1/2}(2-x)\} + 1-(2-x)\right]\\
&= \frac{3x-2}{4},
\end{align*}
where we use $g_{1/2}(2-x) + g_{-1/2}(2-x) = 2(2-x)$ to obtain the last line.

Considering the variable $Y$, only $B=4$ gives a value of $Y$ between 0 and 1, while $B=1,2,3$ lead to the realisation of $Y$ to lie in the interval $(1,2)$. The same calculations as for $X$ give $\Pr(Y\leq y) = y/4$ for $0<y<1$. For $1<y<2$,
\begin{align*}
\Pr(Y\le y) &= \frac{1}{4}\left[1 + 2\Pr\left(2-\frac{1}{Z_P} \le y\right) + \Pr\left(2-\frac{1}{\sqrt{Z_P}}\le y\right)\right] \\
&= \frac{1}{4}\left[1 + 2\{1-(2-y)\} + \{1 - (2-y)^2\}\right] \\
&= \frac{y}{2} - \frac{(2-y)^2}{4}.
\end{align*}

\section*{Declarations}

\subsection*{Ethical Approval} Not applicable

\subsection*{Availability of supporting data} 
Not applicable

\subsection*{Competing Interests} 
No competing interests.

\subsection*{Funding} Not applicable

\subsection*{Authors' contributions} All authors wrote the main manuscript text and C.R. prepared figures 1-5. All authors reviewed the manuscript.

\subsection*{Acknowledgments} Not applicable

\singlespacing
\bibliographystyle{apalike}
\bibliography{sample}

\begin{thebibliography}{}

\bibitem[Beirlant et~al., 2004]{Beirlant2004}
Beirlant, J., Goegebeur, Y., Segers, J., and Teugels, J.~L. (2004).
\newblock {\em {Statistics of Extremes: Theory and Applications}}.
\newblock John Wiley \& Sons Chichester.

\bibitem[Coles et~al., 1999]{ColesHeffernanTawn1999}
Coles, S., Heffernan, J.~E., and Tawn, J.~A. (1999).
\newblock Dependence measures for extreme value analyses.
\newblock {\em Extremes}, 2(4):339--365.

\bibitem[Coles, 2001]{Coles2001}
Coles, S.~G. (2001).
\newblock {\em {An Introduction to Statistical Modeling of Extreme Values}}.
\newblock Springer-Verlag London.

\bibitem[Coles and Tawn, 1994]{Coles1994}
Coles, S.~G. and Tawn, J.~A. (1994).
\newblock Statistical methods for multivariate extremes: An application to
  structural design (with discussion).
\newblock {\em Journal of the Royal Statistical Society: Series C (Applied
  Statistics)}, 43(1):1--49.

\bibitem[Das and Resnick, 2011]{Das2011}
Das, B. and Resnick, S.~I. (2011).
\newblock Conditioning on an extreme component: Model consistency with regular
  variation on cones.
\newblock {\em Bernoulli}, 17(1):226--252.

\bibitem[Davison et~al., 2012]{Davison2012}
Davison, A.~C., Padoan, S.~A., and Ribatet, M. (2012).
\newblock Statistical modeling of spatial extremes.
\newblock {\em Statistical Science}, 27(2):161--186.

\bibitem[Davison and Smith, 1990]{DavisonSmith1990}
Davison, A.~C. and Smith, R.~L. (1990).
\newblock Models for exceedances over high thresholds (with discussion).
\newblock {\em Journal of the Royal Statistical Society: Series B
  (Methodological)}, 52(3):393--425.

\bibitem[Drees and Jan{\ss}en, 2017]{Drees2017}
Drees, H. and Jan{\ss}en, A. (2017).
\newblock Conditional extreme value models: fallacies and pitfalls.
\newblock {\em Extremes}, 20(4):777--805.

\bibitem[Engelke and Hitz, 2020]{Engelke2020}
Engelke, S. and Hitz, A.~S. (2020).
\newblock Graphical models for extremes (with discussion).
\newblock {\em Journal of the Royal Statistical Society: Series B (Statistical
  Methodology)}, 82(4):871--932.

\bibitem[Ewans and Jonathan, 2014]{Ewans2014}
Ewans, K. and Jonathan, P. (2014).
\newblock Evaluating environmental joint extremes for the offshore industry
  using the conditional extremes model.
\newblock {\em Journal of Marine Systems}, 130:124--130.

\bibitem[Gouldby et~al., 2017]{Gouldby2017}
Gouldby, B., Wyncoll, D., Panzeri, M., Franklin, M., Hunt, T., Hames, D.,
  Tozer, N., Hawkes, P., Dornbusch, U., and Pullen, T. (2017).
\newblock Multivariate extreme value modelling of sea conditions around the
  coast of {England}.
\newblock {\em Proceedings of the Institution of Civil Engineers - Maritime
  Engineering}, 170(1):3--20.

\bibitem[Gudendorf and Segers, 2012]{Gudendorf2012}
Gudendorf, G. and Segers, J. (2012).
\newblock Nonparametric estimation of multivariate extreme-value copulas.
\newblock {\em Journal of Statistical Planning and Inference},
  142(12):3073--3085.

\bibitem[Heffernan and Resnick, 2007]{HeffernanResnick2007}
Heffernan, J.~E. and Resnick, S.~I. (2007).
\newblock {Limit laws for random vectors with an extreme component}.
\newblock {\em The Annals of Applied Probability}, 17(2):537 -- 571.

\bibitem[Heffernan and Tawn, 2004]{HT2004}
Heffernan, J.~E. and Tawn, J.~A. (2004).
\newblock A conditional approach for multivariate extreme values (with
  discussion).
\newblock {\em Journal of the Royal Statistical Society: Series B (Statistical
  Methodology)}, 66(3):497--546.

\bibitem[Huser and Wadsworth, 2019]{HuserWadsworth2019}
Huser, R. and Wadsworth, J.~L. (2019).
\newblock Modeling spatial processes with unknown extremal dependence class.
\newblock {\em Journal of the American Statistical Association},
  114(525):434--444.

\bibitem[H{\"u}sler and Reiss, 1989]{Husler1989}
H{\"u}sler, J. and Reiss, R.-D. (1989).
\newblock Maxima of normal random vectors: between independence and complete
  dependence.
\newblock {\em Statistics \& Probability Letters}, 7(4):283--286.

\bibitem[Joe, 1997]{Joe1997}
Joe, H. (1997).
\newblock {\em {Multivariate Models and Multivariate Dependence Concepts}}.
\newblock CRC press.

\bibitem[Keef et~al., 2013a]{Keef2013a}
Keef, C., Papastathopoulos, I., and Tawn, J.~A. (2013a).
\newblock Estimation of the conditional distribution of a multivariate variable
  given that one of its components is large: {Additional} constraints for the
  {Heffernan} and {Tawn} model.
\newblock {\em Journal of Multivariate Analysis}, 115:396--404.

\bibitem[Keef et~al., 2013b]{Keef2013b}
Keef, C., Tawn, J.~A., and Lamb, R. (2013b).
\newblock Estimating the probability of widespread flood events.
\newblock {\em Environmetrics}, 24(1):13--21.

\bibitem[Kiriliouk et~al., 2019]{Kiriliouk2019}
Kiriliouk, A., Rootzén, H., Segers, J., and Wadsworth, J.~L. (2019).
\newblock Peaks over thresholds modeling with multivariate generalized {Pareto}
  distributions.
\newblock {\em Technometrics}, 61(1):123--135.

\bibitem[Nelsen, 1999]{Nelsen1999}
Nelsen, R.~B. (1999).
\newblock {\em {An Introduction to Copulas}}.
\newblock Springer Verlag, New York.

\bibitem[Papastathopoulos et~al., 2017]{Papastathopoulos2017}
Papastathopoulos, I., Strokorb, K., Tawn, J.~A., and Butler, A. (2017).
\newblock Extreme events of {Markov} chains.
\newblock {\em Advances in Applied Probability}, 49(1):134--161.

\bibitem[Papastathopoulos and Tawn, 2016]{Papastathopoulos2016}
Papastathopoulos, I. and Tawn, J.~A. (2016).
\newblock Conditioned limit laws for inverted max-stable processes.
\newblock {\em Journal of Multivariate Analysis}, 150:214--228.

\bibitem[Paulo et~al., 2006]{Paulo2006}
Paulo, M., {van der Voet}, H., Wood, J., Marion, G., and {van Klaveren}, J.
  (2006).
\newblock Analysis of multivariate extreme intakes of food chemicals.
\newblock {\em Food and Chemical Toxicology}, 44(7):994--1005.

\bibitem[Resnick and Zeber, 2014]{Resnick2014}
Resnick, S.~I. and Zeber, D. (2014).
\newblock Transition kernels and the conditional extreme value model.
\newblock {\em Extremes}, 17(2):263--287.

\bibitem[Simpson and Wadsworth, 2021]{Simpson2021}
Simpson, E.~S. and Wadsworth, J.~L. (2021).
\newblock Conditional modelling of spatio-temporal extremes for {Red Sea}
  surface temperatures.
\newblock {\em Spatial Statistics}, 41:100482.

\bibitem[Southworth and Heffernan, 2012]{Southworth2012}
Southworth, H. and Heffernan, J.~E. (2012).
\newblock Extreme value modelling of laboratory safety data from clinical
  studies.
\newblock {\em Pharmaceutical Statistics}, 11(5):361--366.

\bibitem[Tawn, 1990]{Tawn1990}
Tawn, J.~A. (1990).
\newblock Modelling multivariate extreme value distributions.
\newblock {\em Biometrika}, 77(2):245--253.

\bibitem[Tendijck et~al., 2021]{Tendijck2021}
Tendijck, S., Eastoe, E., Tawn, J., Randell, D., and Jonathan, P. (2021).
\newblock Modeling the extremes of bivariate mixture distributions with
  application to oceanographic data.
\newblock {\em Journal of the American Statistical Association}.
\newblock In press.

\bibitem[Wadsworth and Tawn, 2022]{WadsworthTawn2022}
Wadsworth, J. and Tawn, J. (2022).
\newblock Higher-dimensional spatial extremes via single-site conditioning.
\newblock {\em Spatial Statistics}, 51:100677.

\bibitem[Wadsworth et~al., 2017]{Wadsworth2017}
Wadsworth, J.~L., Tawn, J.~A., Davison, A.~C., and Elton, D.~M. (2017).
\newblock Modelling across extremal dependence classes.
\newblock {\em Journal of the Royal Statistical Society: Series B (Statistical
  Methodology)}, 79(1):149--175.

\bibitem[Winter and Tawn, 2017]{WinterTawn2017}
Winter, H.~C. and Tawn, J.~A. (2017).
\newblock kth-order {Markov} extremal models for assessing heatwave risks.
\newblock {\em Extremes}, 20(2):393--415.

\end{thebibliography}

\end{document}